\documentclass[reqno,12pt]{amsart}
\usepackage{latexsym,amssymb,amsfonts, amsmath}
\usepackage{a4wide}
\newtheorem{theorem}{Theorem}[section]
\newtheorem{corollary}[theorem]{Corollary}
\newtheorem{lemma}[theorem]{Lemma}
\newtheorem{proposition}[theorem]{Proposition}
\theoremstyle{definition}

\newtheorem{remark}[theorem]{Remark}

\numberwithin{equation}{section}

\newcommand{\Real}{\mathbb R}
\newcommand{\eps}{\varepsilon}

\newcommand{\ee}{\mathbb{E}}

\newcommand{\nn}{\mathbb{N}}
\newcommand{\rr}{\mathbb{R}}
\newcommand{\pp}{\mathbb{P}}

\newcommand{\zz}{\mathbb{Z}}
\newcommand{\ds}{\displaystyle}

\def\AA{\mathcal A}
\def\BB{\mathcal B}
\def\CC{\mathcal C}

\def\NN{\mathcal N}
\def\PP{\mathcal P}

\def\R{\Real}
\def\N{\nn}

\def\tilde{\widetilde}
\def\hat{\widehat}

\def\dist{{\rm dist}}
\def\const{{\rm const}}
\def\Lip{{\rm Lip}}
\def\hld{{\mathcal H}} %Holder 

\def\ln{\log}
\def\cal{\mathcal}
\def\sm{\smallskip}
\def\med{\medskip}

\def\vep{\varepsilon}
\def\var{\vep}
\def\<{\langle}
\def\>{\rangle}

\def\vphi{\varphi}
\def\derpar#1#2{\ds \frac{\partial #1}{\partial #2}}

\def\beq{\begin{equation}}
\def\eeq{\end{equation}}

\def\beqar{\begin{eqnarray}}
\def\neqar{\end{eqnarray}}

\def\beqa*{\begin{eqnarray*}}
\def\neqa*{\end{eqnarray*}}

\def\bdef{\begin{defn}}
\def\ndef{\end{defn}}

\def\bthm{\begin{theorem}}
\def\nthm{\end{theorem}}

\def\bprop{\begin{prop}}
\def\nprop{\end{prop}}

\def\brmk{\begin{remarks}}
\def\nrmk{\end{remarks}}

\def\bexa{\begin{exa}}
\def\nexa{\end{exa}}

\def\blem{\begin{lem}}
\def\nlem{\end{lem}}

\def\bcor{\begin{cor}}
\def\ncor{\end{cor}}

\def\bexe{\begin{exe}}
\def\nexe{\end{exe}}

\def\bprf{\begin{proof}}
\def\nprf{\end{proof}}

\def\bdes{\begin{description}}
\def\ndes{\end{description}}

\begin{document}
\bibliographystyle{acm}

% ----------------------------------------------------------------
\author{Fran\c cois Bolley}
\address{ENS Lyon, Umpa, 46 all\'ee d'Italie, F-69364 Lyon Cedex 07}
\email{fbolley@umpa.ens-lyon.fr}

\author{Arnaud Guillin}
\address{CEREMADE, Universit\'e Paris Dauphine}
\email{guillin@ceremade.dauphine.fr}

\author{C\'edric Villani}
\address{ENS Lyon, Umpa, 46 all\'ee d'Italie, F-69364 Lyon Cedex 07}
\email{cvillani@umpa.ens-lyon.fr}

\title[]
{Quantitative concentration inequalities for empirical measures on
non-compact spaces}

\subjclass{} \keywords{Transport inequalities, Sanov Theorem}
\thanks{}

\begin{abstract}
	
We establish some quantitative concentration estimates for the empirical 
measure of many independent variables, in transportation distances. 
As an application, we provide some error bounds for particle simulations
in a model mean field problem. The tools include coupling arguments, as
well as regularity and moments estimates for solutions of certain diffusive
partial differential equations.
\end{abstract}

\maketitle
\tableofcontents

\section*{Introduction}
\label{sec-1}

Large stochastic particle systems constitute a popular way to perform
numerical simulations in many contexts, either because they are used in 
some physical model (as in e.g. stellar or granular media) or as an
approximation of a continuous model (as in e.g. vortex simulation for Euler
equation, see \cite[Chapter 5]{MP94} for instance). For such systems one may 
wish to establish concentration estimates
showing that the behavior of the system is sharply stabilized as the
number $N$ of particles goes to infinity. It is natural to search for these
estimates in the setting of large (or moderate) deviations, since one wishes
to make sure that the numerical method has a very small probability to give 
wrong results. From a physical perspective, concentration estimates may be
useful to establish the validity of a continuous approximation such as
a mean-field limit.

When one is interested in the asymptotic behavior of just one, or a few
observables (such as the mean position...), there are efficient methods, based 
for instance on concentration of measure theory. As a good example,
Malrieu~\cite{Mal01} recently applied tools from the fields of Logarithmic Sobolev
inequalities, optimal transportation and concentration of measure,
to prove very neat bounds like
\begin{equation}\label{est1obs}
\sup_{\|\vphi\|_\Lip\leq 1} \pp \left[ \Bigl | 
\frac{1}{N} \sum_{i=1}^{N} \vphi (X_t^i) - \int\vphi\,d\mu_t \Bigr | > \var \right]
\leq 2\, e^{-\lambda N \var^2}. 
\end{equation}

Here $(X_t^i)_{1\leq i\leq N}$ stand for the positions of particles (in phase
space) at time $t$, $\var$ is a given error, $\pp$ stands for the probability,
$\mu_t$ is a probability measure governing the limit behavior of the system,
and $\lambda>0$ is a positive constant depending
on the particular system he is considering (a simple instance of McKean-Vlasov
model used in particular in the modelling of granular media). Moreover,
\[ \|\vphi\|_\Lip := \sup_{x\neq y} \frac{|f(x)-f(y)|}{d(x,y)}, \]
where $d$ is the distance in phase space (say the Euclidean norm $|\cdot|$ in $\R^d$).

This approach can lead to nice bounds, but has the drawback to be limited
to a finite number of observables. Of course, one may apply~\eqref{est1obs}
to many functions $\vphi$, and obtain something like
\begin{equation}\label{estNobs}
\pp \left [ \sum_{k=1}^\infty \frac1{k^2} 
\Bigl | \frac1{N} \sum_{i=1}^N \vphi_k(X_t^i) - \int\vphi_k\,d\mu_t \right |
> \var \Bigr ] \leq C e^{-N\lambda\var^2},
\end{equation}
where $(\vphi_k)_{k\in\N}$ is an arbitrarily chosen dense family in the
set of all $1$-Lipschitz functions converging to~0 at infinity. If we denote by $\delta_x$ 
the Dirac mass at point $x$, and by
\[ \hat{\mu}_t^N := \frac1{N} \sum_{i=1}^N \delta_{X_t^i} \]
the empirical measure associated with the system
(this is a random probability measure), then estimate~\eqref{estNobs} can be
interpreted as a bound on how close $\hat{\mu}_t^N$ is to $\mu_t$. Indeed,
\begin{equation}\label{vweakdist}
d(\mu,\nu) := \sum_{k=1}^\infty \frac1{k^2} \left | 
\int \vphi_k\,d(\mu-\nu) \right |
\end{equation}
defines a distance on probability measures, associated with a topology which
is at least as strong as the weak convergence of measures (convergence against
bounded continuous test functions).
However, this point of view is deceiving: for practical purposes, the distance
$d$ can hardly be estimated, and in any case~\eqref{estNobs} does not
contain more information than~\eqref{est1obs}: it is only useful if one
considers a finite number of observables.

Sanov's large deviation principle~\cite[Theorem 6.2.10]{DZ98} provides 
a more satisfactory tool to estimate the distance between the empirical
measure and its limit. Roughly speaking, it implies, for independent
variables $X_t^i$, an estimate of the form
\[ \pp \left [ \dist (\hat{\mu}^N_t, \mu) \geq \var \right ]
\simeq e^{-N \alpha(\var)}\qquad \text{as $N\to\infty$},\]
where 
\begin{equation}\label{infH}
\alpha(\var) := \inf \Bigl \{ H(\nu|\mu);\: \dist (\nu,\mu) \geq \var \Bigr\}
\end{equation}
and $H$ is the relative $H$ functional:
\[ H(\nu | \mu) = \int \frac{d\nu}{d\mu} \ln \frac{d\nu}{d\mu} \, d\mu \]
(to be interpreted as $+\infty$ if $\nu$ is not absolutely continuous
with respect to $\mu$).
Since $H$ behaves in many ways like a square distance, one can hope that
$\alpha(\var)\geq \const .\, \var^2$. Here ``$\dist$'' may be any distance which is
continuous with respect to the weak topology, a condition which might cause trouble
on a non-compact phase space. 

Yet Sanov's theorem is not the final answer either: it 
is actually asymptotic, and only implies a bound like
\[  \limsup \frac1{N}\log \pp \left [ \dist (\hat{\mu}^N_t, \mu) \geq \var \right ]
\leq -\alpha(\var), \]
which, unlike~\eqref{est1obs}, does not contain any explicit estimate for a given $N$. 
Fortunately, there are known techniques to obtain quantitative upper bounds
for such theorems, see in particular~\cite[Exercise~4.5.5]{DZ98}. Since
these techniques are devised for compact phase spaces, a further truncation will be
necessary to treat more general situations.

In this paper, we shall show how to combine these ideas with recent results
about measure concentration and transportation distances, in order to derive in a
systematic way estimates that are explicit, deal with the empirical measure as a whole, 
apply to non-compact phase spaces, and can be used to study some particle systems
arising in practical problems. Typical estimates will be of the form
\begin{equation}\label{newest}
\pp \left [ \sup_{\|\vphi\|_\Lip\leq 1} \Bigl( \frac1{N} \sum_{i=1}^N
\vphi(X_t^i) - \int\vphi\,d\mu_t \Bigr) > \var \right ] \leq C e^{-\lambda N \var^2}.
\end{equation}
As a price to pay, the constant $C$ in the right-hand side will be much larger than
the one in~\eqref{est1obs}.

Here is a possible application of~\eqref{newest} in a numerical perspective. 
Suppose your system has a limit invariant measure $\mu_\infty=\lim \mu_t$ 
as $t\to\infty$, and you wish to numerically plot its density $f_\infty$. 
For that, you run your particle simulation for a long time $t=T$, and plot, say,
\begin{equation} \label{mutilde}
\tilde{f}_t (x) := \frac1{N} \sum_{i=1}^N \zeta_\alpha \bigl (x-X_t^i\bigr ), 
\end{equation}
where $\zeta_\alpha=\alpha^{-d}\zeta(x/\alpha)$ is a smooth approximation of a Dirac mass
as $\alpha\to 0$ (as usual, $\zeta$ is a nonnegative smooth radial function on $\R^d$ with compact
support and unit integral). With the help of estimates such as~\eqref{newest}, it is often 
possible to compute bounds on, say,
\[ \pp \Bigl [ \|\tilde{f}_T - f_\infty\|_{L^\infty} > \var \Bigr ]\]
in terms of $N$, $\var$, $T$ and $\alpha$. In this way one can ``guarantee''
that all details of the invariant measure are captured by the stochastic system.
While this problem is too general to be treated abstractly, we shall show on some
concrete model examples how to derive such bounds for the same kind of systems that
was considered by Malrieu.

In the next section, we shall explain about our main tools and results;
the rest of the paper will be devoted to the proofs. Some auxiliary
estimates of general interest are postponed in Appendix.

\section{Tools and main results}

\subsection{Wasserstein distances}

To measure distances between probability measures, we shall use transportation distances,
also called {\bf Wasserstein distances}. They can be defined in an abstract Polish
space $X$ as follows: given $p$ in  $[1, +\infty)$,  $d$ a lower semi-continuous 
distance on $X$, and $\mu$ and $\nu$ two Borel probability measures on $X$, the 
Wasserstein distance of order $p$ between $\mu$ and $\nu$ is
$$
 W_p(\mu,\nu):=\inf_{\pi\in \Pi(\mu,\nu)} \left( \iint d(x,y)^p d\pi(x,y) \right)^{1/p} 
$$
where $\pi$ runs over the set $\Pi(\mu,\nu)$ of all joint probability measures on the 
product space $X \times X$ with marginals $\mu$ and $\nu$;
it is easy to check~\cite[Theorem 7.3]{V03} that $W_p$ is a distance on the set $P_p(X)$ of 
Borel probability measures $\mu$ on $X$ such that $\int d(x_0,x)^p\,d\mu(x)<+\infty$.

For this choice of distance, in view of Sanov's theorem, a very natural class of inequalities 
is the family of so-called transportation inequalities, or
{\bf Talagrand inequalities} (see~\cite{Led01} for instance): by definition, given $p \ge 1$ 
and $\lambda>0$, a probability measure $\mu$ on $X$ satisfies $T_p(\lambda)$ if the inequality
$$
W_p(\nu,\mu)\le \sqrt{\frac{2}{\lambda} \, H(\nu \vert \mu)}
$$
holds for any probability measure $\nu$. We shall say that $\mu$ satisfies a $T_p$ inequality
if it satisfies $T_p(\lambda)$ for some $\lambda >0$. By Jensen's inequality, these
inequalities become stronger as $p$ becomes larger; so the weakest of all is $T_1$.
Some variants introduced in~\cite{BV04} will also be considered.

Of course $T_p$ is not a very explicit condition, and a priori it is not clear how to
check that a given probability measure satisfies it. It has been 
proven~\cite{BG99,DGW02,BV04} that $T_1$ is \emph{equivalent} to the existence of a 
square-exponential moment: in other words, a reference measure $\mu$ satisfies $T_1$ if and 
only if there is $\alpha>0$ such that
\[ \int e^{\alpha d(x,y)^2} \, d\mu(x) < +\infty \]
for some (and thus any) $y\in X$. If that condition is satisfied, then one can find explicitly
some $\lambda$ such that $T_1(\lambda)$ holds true: see for instance~\cite{BV04}.

This criterion makes $T_1$ a rather convenient inequality to use. Another popular
inequality is $T_2$, which appears naturally in many situations where a lot of structure
is available, and which has good tensorization properties in many dimensions.
Up to now, $T_2$ inequalities have not been so well characterized: it is known that they
are implied by a Logarithmic Sobolev inequality~\cite{OV00,BGL01,Wan03}, and that they
imply a Poincar\'e, or spectral gap, inequality~\cite{OV00,BGL01}. See~\cite{CG04} for an attempt 
to a criterion for $T_2$.
In any case, contrary to the case $p=1$, there is no hope to obtain $T_2$ inequalities 
from just integrability or decay estimates. 

In this paper, we shall mainly focus on the case $p=1$, which is much more flexible.

\subsection{Metric entropy}

When $X$ is a compact space, the minimum number $m(X,r)$ of balls of
radius $r$ needed to cover $X$ is called the {\bf metric entropy} of $X$. 
This quantity plays an important role in quantitative variants of Sanov's 
Theorem~\cite[Exercise 4.5.5]{DZ98}. In the present paper, to fix ideas we shall
always be working in the particular Euclidean space $\R^d$, which of course is not compact;
and we shall reduce to the compact case by truncating everything to balls of finite radius $R$.
This particular choice will influence the results through the function $m(\PP_p(B_R),r)$, where
$B_R$ is the ball of radius $R$ centered at some point, say the origin, and $\PP_p(B_R)$
is the space of probability measures on $B_R$, metrized by $W_p$.

\subsection{Sanov-type theorems}

The core of our estimates is based on variants of Sanov's Theorem, all dealing with
\emph{independent} random variables. Let $\mu$ be a given probability measure on
$\R^d$, and let $(X^i)_{i=1,...,N}$ be a sample of independent variables, all distributed
according to $\mu$; let also 
\[ \hat{\mu}^N := \frac1{N} \sum_{i=1}^N \delta_{X^i}\]
be the associated empirical measure. In our first main result we assume a $T_p$
inequality for the measure $\mu$, and deduce from that an upper bound in $W_p$ distance:

\bthm\label{thmconc-Tp2}
Let $p \in [1,2]$ and let $\mu$ be a probability measure on $\rr^d$ satisfying a 
$T_p(\lambda)$ inequality. Then, for any $d'>d$ and $\lambda' < \lambda$, there exists 
some constant $N_0$, depending on $\lambda',  d'$ and some square-exponential moment of $\mu$, such that
for any $\eps >0$ and $N \geq N_0 \max(\eps^{-(d'+2)},1)$,
\begin{equation}\label{conc-Tp2} 
\pp\left[ W_p(\mu,\hat{\mu}^N)>\eps\right]\le e^{- \gamma_p \, \frac{\lambda'}{2} \, N \, \eps^2},
\end{equation}
where
\[ \gamma_p = \begin{cases} 1 \text{ if $1\leq p < 2$ } \\
3 - 2 \, \sqrt{2} \text{ if $p=2$ }. \end{cases} \]
\nthm

Compared to Sanov's Theorem, this result is more restrictive in the sense that it requires some
extra assumptions on the reference measure $\mu$, but under these hypotheses we are able to
replace a result which was only asymptotic by a pointwise upper bound on the error probability, 
together with a lower bound on the required size of the sample.

In view of the Kantorovich-Rubinstein duality formula
\begin{equation}\label{KR}
W_1(\mu,\nu)=\sup\left\{\int f \,d(\mu-\nu);\:\: \| f\|_{\Lip}\le 1\right\},
\end{equation}
Theorem~\ref{thmconc-Tp2} implies concentration inequalities such as
$$
\pp\left[\sup_{f;\:\| f\|_{\Lip}\le 1}\Bigl(\frac{1}{N}\sum_{k=1}^{N} f(X_i) -
\int f \, d\mu\Bigr) >\vep \right]\le e^ {-\frac{\lambda'}{2}N\vep^ 2}
$$
for $\lambda'<\lambda$, and $N$ sufficiently large, under the assumption that $\mu$
satisfies a $T_1$ inequality, or equivalently admits a finite square-exponential moment.
Those types of inequalities are of interest in non-parametric statistics and choice 
models~\cite{Mas03}. 

\begin{remark} The sole inequality $T_1(\lambda)$ implies that for all 1-Lipschitz function $f$,
$$\pp\left[\frac{1}{N}\sum_{k=1}^{N} f(X_i) -\int f \, d\mu >\vep \right ]\le 
e^ {-{\lambda\over 2}N\vep^ 2},$$
and it is easy to see that the coefficient $\lambda$ in this inequality is the best
possible. While the quantity controlled in Theorem~\ref{thmconc-Tp2} is much stronger,
the estimate is weakened only in that $\lambda$ is replaced by some $\lambda'>\lambda$
(arbitrarily close to $\lambda$) and that $N$ has to be large enough.
In fact, a variant of the proof below would yield estimates such as
$$\pp\left[ W_p(\mu,\hat{\mu}^N)>\eps\right]\le C(\eps)\,
e^{- \gamma \, \frac{\lambda'}{2} \, N \, \eps^2},
$$
where now there is no restriction on $N$, but $C(\eps)$ is a larger constant, explicitly
computable from the proof.
\end{remark}

\begin{remark} As pointed out to us by M. Ledoux, there is another way to concentration estimates on the 
empirical measure when $d=p=1$. Indeed, in this specific case,
\[
W_1(\hat{\mu}^N, \mu) = \Big\Vert \frac{1}{N} \sum_{i=1}^{N} H(\cdot - X_i) - F \Big\Vert_{L^1(\rr)}
\]
where $H = {\bf 1}_{[0,+\infty)}$ stands for the Heaviside function on $\rr$ and $F$ denotes the
repartition function of $\mu$, so that
\[
\pp \Big[ W_1(\hat{\mu}^N, \mu) \geq \eps \Big]  = 
\pp\left[ \Big\Vert \frac{1}{N} \sum_{i=1}^{N} F_i \Big\Vert_{L^1} >\eps \right]
\]
where
\[
F_i := H(\cdot - X_i) - F \qquad (1 \leq i \leq N)
\]
are centered $L^1(\rr)$-valued independent identically distributed random variables. But, according to 
\cite[Exercise 3.8.14]{AZ80}, a centered  $L^1(\rr)$-valued random variable $Y$ satisfies a Central 
Limit Theorem if and only if
\[
\int_{\rr} \left( \ee[Y^2(t)] \right)^{1/2} \, dt < + \infty,
\]
a condition which for the random variables $F_i$'s can be written
\begin{equation}\label{condCLT}
\int_{\rr} \sqrt{F(t) (1-F(t))} \, dt < +\infty.
\end{equation}
Condition \eqref{condCLT} in turn holds true as soon as (for instance) 
$\int_{\rr} \vert x \vert^{2+\delta} \, d\mu(x) $ is finite for some positive $\delta$. 
Then we may apply a quantitative version of the Central Limit Theorem for random varaiables
in the Banach space $L^1(\rr)$. See \cite{GZ86} and \cite{LT91} for related works.

\end{remark}

\begin{remark} Theorem~\ref{thmconc-Tp2} applies if $N$ is at least as
large as $\var^{-r}$ for some $r>d+2$; we do not know whether $d+2$
here is optimal.
\end{remark}

For the applications that we shall treat, in which the tails of the probability
distributions will be decaying very fast, Theorem~\ref{thmconc-Tp2} will be sufficient.
However, it is worthwile pointing out that the technique works under much broader assumptions:
weaker estimates can be proven for probability measures that do not decay fast enough to admit 
finite square-exponential moments. Here below are some such results using only polynomial moment 
estimates:

\bthm\label{thmconc-Mq}
Let $q\geq 1$ and let $\mu$ be a probability measure on $\rr^d$ such that 
\[\int_{\rr^d} \vert x \vert^q \, d\mu(x) < +\infty. \] 
Then

\sm
(i) For any $p \in [1,q/2)$, $\delta \in (0,q/p -2)$ and $d'>d$, there exists a constant $N_0$ 
such that
$$ 
\pp\left[ W_p(\mu,\hat{\mu}^N)>\eps\right]\le \eps^{-q} N^{-\frac{q}{2p} + \frac{\delta}2} 
$$
for any $\eps>0$ and $N \geq N_0\, \max(\eps^{-q \frac{2p+d'}{q-p}},\eps^{d'-d})$;
\sm

(ii)  For any $p \in [q/2,q)$, $\delta \in (0,q/p -1)$ and $d'>d$ there exists a constant $N_0$ 
such that
$$ 
\pp\left[ W_p(\mu,\hat{\mu}^N)>\eps\right]\le \eps^{-q} N^{1 - \frac{q}p + \delta} 
$$
for any $\eps>0$ and $N \geq N_0 \max(\eps^{-q \frac{2p+d'}{q-p}},\eps^{d'-d})$.

\nthm

%\begin{remark}
%Again, variants are possible: for instance, estimate (ii) could be stated, up to some 
%multiplicative constant, for any $\delta>0$, but then, for a given $\var$ the right-hand side 
%would a priori go to~0 as $N\to\infty$ only if $\delta<(q/p)-1$.
%Moreover, by setting $p=1$ and letting $q$ go to $+\infty$, we recover the condition 
%$N \geq N_0\,\eps^{-(2+d')}$ in Theorem~\ref{thmconc-Tp2}. Finally, by setting  $q = 2\, p$, in 
%particular for $p=1$ and $q=2$, we can bound the probability of error by $(N \, \eps^2)^{-1}$.
%\end{remark}
%
%{\bf ndFB: remarque a laisser sous cette forme ?}

Here are also some variants under alternative ``regularity'' assumptions:

\bthm \label{thmconc-var}
\begin{enumerate}
\item[(i)] Let $p \geq 1$; assume that ${\mathcal E}_{\alpha} := \int e^ {\alpha \vert x \vert} d\mu$ 
is finite for some $\alpha > 0$. Then, for all $d' > d$, there exist
some constants $K$ and $N_0$, depending only on $d$, $\alpha$ and ${\mathcal E}_{\alpha}$, such that
\[ 
\pp\left[ W_p(\mu,\hat{\mu}^N)>\eps\right]\le e^{-  K \, N^{1/p} \, \min(\eps,\eps^2)}
\]
for any $\eps >0$ and $N \geq N_0 \max(\eps^{-(2 p  + d')},1)$.
\sm

%for all $d' > d$, there exists a constant $K$, depending
%only on $d$, $\alpha$ and ${\cal E}_\alpha$, such that
%\[ 
%\pp\left[ W_p(\mu,\hat{\mu}^N)>\eps\right]\le \exp \left( \frac{\eps^{-d'}}{K} -  
%K \, N \, \min(\eps,\eps^2)^p \right)
%\]
%for any $\eps >0$ and $N \ge 1$.

\item[(ii)] Suppose that $\mu$ satifies $T_1$ and a Poincar\'e inequality, then for all $a<2$ there exists 
some constants $K$ and $N_0$ such that
\begin{equation} 
\pp\left[ W_2(\mu,\hat{\mu}^N)>\eps\right]\le e^{- K \, N \, \min(\eps^2,\eps^a)}
\end{equation}
for any $\eps >0$ and $N \geq N_0 \max(\eps^{-(4  + d')},1)$.
\sm

\item[(iii)]  Let $p>2$ and let $\mu$ be a probability measure on $\rr^d$ satisfying $T_p(\lambda)$.
Then for all $\lambda' < \lambda$ and $d'>d$ there exists some constant $N_0$, depending 
on $\mu$ only through $\lambda$ and some square-exponential moment, such that
\begin{equation}\label{conc-Tp} 
\pp\left[ W_p(\mu,\hat{\mu}^N)>\eps\right]\le
\min \left( e^{-\frac{\lambda'}{2} N\eps^2} + e^{-(N\eps^{d'+2})^{2/d'}} , \, 2 \, 
e^{-\frac{\lambda'}{4} N^{2/p} \eps^2} \right)
\end{equation}
for any $\eps >0$ and $N \geq N_0\,\max(\eps^{-(d'+2)},1)$.
\end{enumerate}
\nthm

\subsection{Interacting systems of particles}

We now consider a system of $N$ interacting particles whose time-evolution is
governed by the system of coupled stochastic differential equations
\begin{equation}\label{edsXi}
dX_t^i = \sqrt{2} \, dB_t^i - \nabla V(X_t^i) dt - \frac{1}{N} \sum_{j=1}^{N} 
\nabla W(X_t^i-X_t^j) dt,
\qquad i=1, \dots, N.
\end{equation}
Here $X_t^i$ is the position at time $t$ of particule number $i$,
the $B^i$'s are $N$ independent Brownian motions, and $V$ and $W$ are smooth potentials,
sufficiently nice that~\eqref{edsXi} can be solved globally in time.
We shall always assume that $W$ (which can be interpreted as an interaction potential)
is a symmetric function, that is $W(-z) = W(z)$ for all $z\in \rr^d$.

Equation~\eqref{edsXi} is a particularly simple instance of coupled system;
in the case when $V$ is quadratic and $W$ has cubic growth, it was used
as a simple mean-field kinetic model for granular media (see e.g.~\cite{Mal01}).
While many of our results could be extended to more general systems,
that particular one will be quite enough for our exposition.

To this system of particles is naturally associated the empirical measure, defined for each 
time $t \geq 0$ by
\begin{equation}\label{empX_t}
\hat{\mu}_t^N \, := \, \sum_{i=1}^{N} \delta_{X_t^i}.
\end{equation}
Under suitable assumptions on the potentials $V$ and $W$,
it is a classical result that, if the initial positions of the particle system
are distributed chaotically (for instance, if they are identically distributed,
independent random variables), then the empirical measure $\hat{\mu}_t^N$ converges
as $N\to\infty$ to a solution of the nonlinear partial differential equation
\begin{equation}\label{edp}
\frac{\partial \mu_t}{\partial t} = \Delta \mu_t + 
\nabla\cdot\Bigl (\mu_t \nabla\bigl(V+W*\mu_t\bigr)\Bigr),
\end{equation}
where $\nabla\cdot$ stands for the divergence operator. Equation~\ref{edp}
is a simple instance of McKean-Vlasov equation. This convergence result is part of 
the by now well-developed theory of propagation of chaos, and was studied
by Sznitman for pedagogical reasons~\cite{Szn91}, in the case of potentials 
that grow at most quadratically at infinity.
Later, Benachour, Roynette, Talay and Vallois~\cite{BRTV98,BRV98} considered
the case where the interaction potential grows faster than quadratically.
As far as the limit equation~\eqref{edp} is concerned, a discussion of its use
in the modelling of granular media in kinetic theory was performed by Benedetto, Caglioti, Carrillo
and Pulvirenti~\cite{BCCP98,BCP97}, while the asymptotic behavior in large time
was studied by Carrillo, McCann and Villani~\cite{CMV03,CMV04} with the help of 
Wasserstein distances and entropy inequality methods. Then Malrieu~\cite{Mal01}
presented a detailed study of both limits $t\to\infty$ and $N\to\infty$ by probabilistic
methods, and established estimates of the type of~\eqref{est1obs} under adequate convexity 
assumptions on $V$ and $W$ (see also~\cite[Problem~15]{V03}).
\medskip

As announced before, we shall now give some estimates on the convergence at the level 
of the law itself. To fix ideas, we assume that $V$ and $W$ have locally bounded
Hessian matrices satisfying
\begin{equation}\label{D2V,D2W,V}
\left\{ \begin{array}{cl}
        \text{(i)} & D^2V(x) \geq \beta I, \quad \gamma I \leq D^2W(x) \leq \gamma' I,  
	\qquad \forall x \in \rr^d,\\ \\
        \text{(ii)} &  |\nabla V(x)| = O(e^{a \vert x \vert^2}) \; \text{\quad for any} \, a >0.  
                          \end{array}
                  \right.
\end{equation}

Under these assumptions, we shall derive the following bounds.

\bthm\label{thmconc-eds}
Let $\mu_0$ be a probability measure on $\rr^d$, admitting a finite square-exponential
moment:
\[ \exists\alpha_0>0;\qquad M_{\alpha_0} := \int e^{\alpha_0 |x|^2}\,d\mu_0(x) <+\infty.\]
Let $(X_0^i)_{1 \leq i \leq N}$ be $N$ independent random variables with common 
law $\mu_0$. Let $(X_t^i)$ be the solution of~\eqref{edsXi} with initial value
$(X_0^1,\ldots X_0^N)$, where $V$ and $W$ are assumed to satisfy~\eqref{D2V,D2W,V}; 
and let $\mu_t$ be the solution of~\eqref{edp} with initial value $\mu_0$. Let
also $\hat{\mu}^N_t$ be the empirical measure associated with the $(X_t^i)_{1\leq i\leq N}$. 
Then, for all $T\geq 0$, there exists some constant $K=K(T)$  such that,
for any $d'>d$, there exists some constants $N_0$ and $C$ such that for all $\var>0$
$$
N \geq N_0\,\max(\varepsilon^{-(d'+2)},1)\Longrightarrow\quad
\pp \, \left[\sup_{0 \leq t \leq T} W_1(\hat{\mu}_t^N,\mu_t) > \varepsilon\right] \leq 
C (1+T\var^{-2}) \exp \left( - K \, N\, \vep^2 \right).
$$
\nthm

Note that in the above theorem we have proven not only that for all $t$, the
empirical measure is close to the limit measure, but also that the probability of
observing any significant deviation during a whole time period $[0,T]$ is small.

The fact that $\hat{\mu}^N_t$ is very close to the deterministic
measure $\mu_t$ implies the propagation of chaos: two particles drawn
from the system behave independently of each other as $N\to\infty$
(see Sznitman~\cite{Szn91} for more details).
But we can also directly study correlations between particles and
find more precise estimates: for that purpose it is convenient to consider
the empirical measure on \emph{pairs} of particles, defined as
\[ \hat{\mu}^{N,2}_t := \frac1{N(N-1)} \sum_{i\neq j} \delta_{(X_t^i,X_t^j)}. \]
By a simple adaptation of the computations appearing in the proof
of Theorem~\ref{thmconc-eds}, one can prove

\bthm \label{thmpropchaos}
With the same notation and assumptions as in Theorem~\ref{thmconc-eds}, for all $T \geq 0$ and $d'>d$,
there exists some constants $K>0$ and $N_0$ such that for all $\eps >0$
\[ N \geq N_0 \max(\varepsilon^{-(d'+2)},1) \Longrightarrow\quad
\pp \, \left[ W_1(\hat{\mu}_t^{N,2},\mu_t\otimes\mu_t) > \varepsilon\right] 
\leq \exp \left( - K \, N\, \vep^2 \right). \]
\nthm

(Here $W_1$ stands for the Wasserstein distance or order $1$ on $P_1(\rr^d \times \rr^d)$.)
Of course, one may similarly consider the problem of drawing $k$ particles with $k\geq 2$.
\med

Theorems~\ref{thmconc-eds} and~\ref{thmpropchaos} 
use Theorem~\ref{thmconc-Tp2} as a crucial ingredient,
which is why a strong integrability assumption is imposed on $\mu_0$.
Note however that, under stronger assumptions on the behaviour at
infinity of $V$ or $W$, as the existence of some $\beta \in \rr$, $B, \, \eps >0$ such as
\[
D^2 V(x) \geq (B \vert x \vert^{\eps} + \beta ) I, \qquad, \forall x \in \rr^d,
\]
it can be proven that any square exponential moment for $\mu_t$ becomes 
instantaneously finite for $t >0$. Note also that, by using Theorem~\ref{thmconc-Mq}, 
one can obtain weaker but still relevant results of concentration of the empirical 
measure under just polynomial moment assumptions on $\mu_0$, provided
that $\nabla V$ does not grow too fast at infinity.
To limit the size of this paper, we shall not go further into
such considerations.

\subsection{Uniform in time estimates}

In the ``uniformly convex case'' when $\beta >0, \beta+ 2 \gamma>0$, 
%or when $\beta + \gamma >0$ and
%the center of mass $\ds \int_{\rr^d} x \, d\mu_t(x)$ is invariant by the evolution along the equation 
%(which is true if $V \equiv 0$ for instance), 
it can be proven~\cite{Mal01,CMV03,CMV04}
that $\mu_t$ converges exponentially fast, as $t\to\infty$, to some
equilibrium measure $\mu_\infty$. In that case, it is natural to expect that the 
empirical measure is a good approximation of $\mu_\infty$ as $N\to\infty$ and $t\to\infty$,
uniformly in time. This is what we shall indeed prove:

%, focusing on the case when the
%convexity of $V$ is strong enough to balance a lack of convexity of $W$; in another case when
%$V \equiv 0$, or more generally when the center of mass is kept fixed, one may hope to obtain
%such results, but we did not go further in that direction. More precisely:

\bthm\label{thmconc-eds-unif}
With the same notation and assumptions as in Theorem~\ref{thmconc-eds},
suppose that $\beta> 0, \beta + 2 \gamma >0$. Then there exists some constant $K>0$ such that
for any $d'>d$, there exists some constants $C$ and $N_0$ such that for all $\var>0$
$$
N \geq N_0 \max(\varepsilon^{-(d'+2)},1) \Longrightarrow\quad
\sup_{t\geq 0} \pp \, [ W_1(\hat{\mu}_t^N,\mu_t) > \varepsilon] \leq 
C (1+\var^{-2})\; \exp \left( - K \, N\, \vep^2 \right) 
$$
As a consequence, there are constants $T_0$, $\var_0$ (depending on the initial datum) and $K'=K/4$
such that, under the same conditions on $N$ and $\var$,
\[ \sup_{t\geq T_0\log (\var_0/\var)} \pp \, [ W_1(\hat{\mu}_t^N,\mu_\infty) > \varepsilon] \leq 
C(1+\var^{-2})\; \exp \left( - K' \, N\, \vep^2 \right). \] 
\nthm

\begin{remark} In view of the results in~\cite{CMV03}, it is natural to
expect that a similar conclusion holds true when $V=0$ and $W$ is convex 
enough. Propositions \ref{propT_1edp} and \ref{propasympt} below extend to that case, but
it seems trickier to adapt the proof of Proposition \ref{propasympt}.

\end{remark}

We conclude with an application to the numerical reconstruction of the invariant measure.

\bthm\label{thmreconstr}
With the same notation and assumptions as in Theorem~\ref{thmconc-eds-unif}, consider
the mollified empirical measure~\eqref{mutilde}. Then one can choose $\alpha=O(\var)$
in such a way that
\begin{multline*}
 N \geq N_0 \max(\varepsilon^{-(d'+2)},1) \Longrightarrow\quad \sup_{t\geq T_0\log(\var_0/\var)} \pp \, 
\Bigl[ \|\tilde{f}_t - f_\infty\|_{L^\infty} > \varepsilon\Bigr] \\
\leq C(1+\var^{-(2d+4)})\; \exp \left( - K' \, N\, \vep^{2d+4} \right). 
\end{multline*}
\nthm

These results are effective: all the constants therein can be estimated explicitly
in terms of the data.

\subsection{Strategy and plan}

The strategy is rather systematic. First, we shall establish Sanov-type bounds for
independent variables in $\R^d$ (not depending on time), resulting in concentration results 
such as Theorems~\ref{thmconc-Tp2} to~\ref{thmconc-var}. This will be
achieved along the ideas in~\cite[Exercices 4.5.5 and~6.2.19]{DZ98} (see also
\cite[Section 5]{Sch96}),
by first truncating to a compact ball, and then covering the set of probability
measures on this ball by a finite number of small balls (in the space of probability 
measures); the most tricky part will actually lie in the optimization of parameters.

With such results in hand, we will start the study of the particle system by
introducing the nonlinear partial differential equation~\eqref{edp}.
For this equation, the Cauchy problem can be solved in a satisfactory way, in 
particular existence and uniqueness of a solution, which for $t>0$ is reasonably
smooth, can be shown under various assumptions on $V$ and $W$ (see e.g.~\cite{CMV03,CMV04}).
Other regularity estimates such as the decay at infinity, or the smoothness in time,
can be established; also the convergence to equilibrium in large time can sometimes
be proven.

Next, following the presentation by Sznitman~\cite{Szn91}, we introduce a family of 
independent processes $(Y_t^i)_{1\leq i\leq N}$, governed by the stochastic differential equation
\begin{equation}\label{nlsde}
\left\{ \begin{array}{rcl}
             dY^i_t & = & \sqrt{2}\,dB^i_t - \nabla V(Y^i_t)\,dt - \nabla W\ast\mu_t(Y^i_t)\,dt,\\
              Y^i_0 & = & X^i_0 .  
                          \end{array}
                  \right.
\end{equation}
As a consequence of It\^o's formula, the law $\nu_t$ of each $Y_t^i$ is a solution of the linear
partial differential equation
\[ \derpar{\nu_t}{t} = \Delta\nu_t + \nabla\cdot 
\Bigl(\nabla \bigl (V+W\ast\mu_t\bigr) \nu_t\Bigr), \qquad \nu_0=\mu_0. \]
But this linear equation is also solved by $\mu_t$, and a uniqueness theorem implies that
actually $\nu_t=\mu_t$, for all $t\geq 0$. See~\cite{BRTV98,BRV98} for related questions on 
the stochastic differential equation~\eqref{nlsde}. 

For each given $t$, the independence of the variables $Y_t^i$ and the good decay of $\mu_t$
will imply a strong concentration of the empirical measure
\[ \hat{\nu}^N_t := \frac1{N} \sum_{i=1}^N \delta_{Y_t^i}. \]
To go further, we shall establish a more precise information, such as a control on
\[ \pp \, \left[\sup_{0 \leq t \leq T} W_1(\hat{\nu}^N_t,\mu_t) > \varepsilon\right]. \]
Such bounds will be obtained by combining the estimate of concentration at fixed time $t$
with some estimates of regularity of $\hat{\nu}^N_t$ (and $\mu_t$) in $t$, obtained
via basic tools of stochastic differential calculus (in particular Doob's inequality).

Finally, we can show by a Gronwall-type argument that the control of the distance
of $\hat{\mu}^N_t$ to $\mu_t$ reduces to the control of the distance of $\hat{\nu}^N_t$ to
$\mu_t$: for instance,
\begin{equation}\label{W1coupl}
\pp \, \left[\sup_{0 \leq t \leq T} W_1(\hat{\mu}_t^N,\mu_t) > \varepsilon\right] \leq
\pp \, \left[\sup_{0 \leq t \leq T} W_1(\hat{\nu}_t^N,\mu_t) > C \varepsilon\right]
\end{equation}
for some constant $C$.
We shall also show how a variant of this computation provides estimates of the
type of those in Theorem~\ref{thmconc-eds-unif}, and how to get data reconstruction
estimates as in Theorem~\ref{thmreconstr}.

\subsection{Remarks and further developments}

The results in this paper confirm what seems to be a rather general 
rule about Wasserstein distances:
results in distance $W_1$ are very robust and can be used in rather hard problems,
with no particular structure; on the contrary, results in distance $W_2$ are
stronger, but usually require much more structure and/or assumptions. For instance,
in the study of the equation~\eqref{edp}, the distance $W_2$ works beautifully,
and this might be explained by the fact that~\eqref{edp} has the structure of a gradient 
flow with respect to the $W_2$ distance~\cite{CMV03,CMV04}. In the problem
considered by Malrieu~\cite{Mal01}, $W_2$ is also well-adapted, but leads him to impose
strong assumptions on the initial datum $\mu_0$, such as the existence of a
Logarithmic Sobolev inequality for $\mu_0$, considered as a reference measure. As a general rule,
in a context of geometric inequalities with more or less subtle isoperimetric content,
related to Brenier's transportation mapping theorem, $W_2$ is also the most natural distance 
to use~\cite{V03}. On the contrary, here we are considering quite a rough problem
(concentration for the law of a random probability measure, driven by a stochastic
differential equation with coupling) and we wish to impose only natural
integrability conditions; then the distance $W_1$ is much more convenient.

Further developments could be considered. For instance, one may desire to prove some
deviation inequalities for dependent sequences, say Markov chains, as both Sanov's 
theorem and transportation inequality can be established under appropriate ergodicity 
and integrability conditions. 

Considering again the problem of the particle system, in a numerical context, one may 
wish to take into account the numerical errors associated with
the time-discretization of the dynamics (say an implicit Euler scheme). 
For concentration estimates in one observable, a beautiful study of these issues was
performed by Malrieu~\cite{Mal03}. For concentration estimates on the whole
empirical measure, to our knowledge the study remains to be done.
Also errors due to the boundedness of the phase space actually used
in the simulation might be taken into account, etc.

At a more technical level, it would be desirable to relax the
assumption of boundedness of $D^2W$ in Theorem~\ref{thmconc-eds},
so as to allow for instance the interesting case of cubic
interaction. This is much more technical and will be considered
in a separate work.

Another issue of interest would be to consider concentration of the empirical
measure \emph{on path space}, i.e. 
\[ \hat{\mu}^N_{[0,T]} := \frac1N \sum_{i=1}^N \delta_{(X_t^i)_{0\leq t\leq T}}, \]
where $T$ is a fixed time length. Here $\hat{\mu}^N_{[0,T]}$ is a random measure
on $C([0,T];\R^d)$ and we would like to show that it is close to the law of
the trajectories of the nonlinear stochastic differential equation
\begin{equation}\label{nlSDE} 
dY_t = \sqrt{2}\,dB_t - \nabla V(Y_t)\,dt - (\nabla W*\mu_t)(Y_t)\,dt, 
\end{equation}
where the initial datum $Y_0$ is drawn randomly according to $\mu_0$.
This will imply a quantitative information on the whole trajectory of a given particle 
in the system. 

When one wishes to adapt the general method to this question, a problem immediately
occurs: not only is $C([0,T];\R^d)$ not compact, but also
balls with finite radius in this space are not compact either (of course, this is
true even if the phase space of particles is compact). One may remedy to this problem
by embedding $C([0,T];B_R)$ into a space such as $L^2([0,T];B_R)$, equipped with
the weak topology; but we do not know of any ``natural'' metric on that space.
There is (at least) another way out: we know from classical stochastic processes theory that
integral trajectories of differential equations driven by white noise are typically
H\"older-$\alpha$ for any $\alpha<1/2$. This suggests a natural strategy:
choose any fixed $\alpha\in (0,1/2)$ and work in the space $\hld^\alpha([0,T];\R^d)$,
equipped with the norm
\[ \|w\|_{\hld^\alpha} := \sup_{0\leq t\leq T} |w(t)| +
\sup_{s\neq t} \frac{|w(t)-w(s)|}{|t-s|^\alpha}. \]
For any $R>0$, the ball of radius $R$ and center~0 (the zero function) in 
$\hld^\alpha$ is compact, and one may estimate its metric entropy. Then one can hope
to perform all estimates by using the norm $\hld^\alpha$; for instance,
establish a bound on, say, a square-exponential moment on the law of $Y_t$:
\[ \ee \exp\left (\beta \|(Y_t)_{0\leq t\leq T}\|_{\hld^\alpha}^2\right ) < +\infty. \]
Again, to avoid expanding the size of the present paper too much, 
these issues will be addressed separately.

\section{The case of independent variables}

In this section we consider the case where we are given $N$ \emph{independent}
variables $X^i\in\R^d$, distributed according to a certain law $\mu$. There is
no time dependence at this stage. We shall first examine the case when the
law $\mu$ has very fast decay (Theorem~\ref{thmconc-Tp2}), then variants in which
it decays in a slower way (Theorem~\ref{thmconc-Mq} and~\ref{thmconc-var}).

\subsection{Proof of Theorem \ref{thmconc-Tp2}}
\label{sectionpreuve}

The proof splits into three steps: (1) Truncation to a compact ball $B_R$ of radius $R$,
(2) covering of $\PP(B_R)$ by small balls of radius $r$ and Sanov's argument, and
(3) optimization of the parameters.
\bigskip

{\bf Step~1: Truncation.} Let $R>0$, to be chosen later on, and let $B_R$ stand for the ball
of radius $R$ and center 0 (say) in $\rr^d$. Let ${\bf 1}_{B_R}$ stand for the indicator 
function of $B_R$. We truncate $\mu$ into a probability measure $\mu_R$ on the ball $B_R$:
$$
\mu_R = \frac{{\bf 1}_{B_R}\,\mu}{\mu[B_R]}.
$$

We wish to bound the quantity $\pp\left[ W_p(\hat{\mu}^N,\mu)>\eps\right]$
in terms of $\mu_R$ and the associated empirical measure.
For this purpose, consider independent variables $(X^k)_{1\leq k\leq N}$ drawn
according to $\mu$, and $(Y^k)_{1\leq k\leq N}$ drawn according to $\mu_R$,
independent of each other; then define
$$X^k_R:=\left\{\begin{array}{ll}
X^k & \text{if } |X^k|\le R\\
Y^k & \text{if } |X^k|> R.\\
\end{array}\right.
$$

Since $X^1$ and $X^1_R$ are distributed according to $\mu$ and $\mu_R$ respectively,
we have, by definition of Wasserstein distance,
\begin{multline*}
W_p^p(\mu,\mu_R) \le \ee|X^1-X^1_R|^p = 
\ee\Bigl( |X^1-Y^1|^p {\bf 1}_{|X^1|>R}\Bigr) 
\leq 2^p \ee \bigl(|X^1|^p {\bf 1}_{|X^1|>R}\bigr) \\
= 2^p \int_{\{\vert x \vert > R\}} \vert x \vert^p \, d\mu(x).
\end{multline*}
But $\mu$ satisfies a $T_p(\lambda)$ inequality for some $p\geq 1$, hence
a fortiori a $T_1(\lambda)$ inequality, so 
\[ E_\alpha:= \int_{\rr^d} e^{\alpha|x|^2}\,d\mu(x) <+\infty \]
for some $\alpha>0$ (any $\alpha<\lambda/2$ would do). If $R$ is large enough
(say, $R\geq \sqrt{p/(2\alpha)}$), then the function
$\ds r \longmapsto \frac{r^p}{e^{\alpha r^2}}$
is nonincreasing for $r\geq R$, and then
\[ W_p^p(\mu,\mu_R) \le 2^p \left (\frac{R^p}{e^{\alpha R^2}} \right )
\int_{\{\vert x \vert > R\}} e^{\alpha |x|^2} \, d\mu(x).\]
We conclude that
\begin{equation}\label{wp-trunc-ineq}
W_p^p(\mu,\mu_R) \le 2^p E_\alpha R^p e^{-\alpha R^2} \qquad
(\alpha<\lambda/2, \; R\ge \sqrt{p/2\alpha}).
\end{equation}

On the other hand, the empirical measures
$$
\hat{\mu}^N:=\frac{1}{N}\sum_{k=1}^N \delta_{X^k},   \qquad \qquad 
\hat{\mu}^N_R:=\frac{1}{N}\sum_{k=1}^N \delta_{X^k_R}
$$
satisfy
$$
W_p^p(\hat{\mu}^N_R,\hat{\mu}^N) \le \frac{1}{N}\sum_{k=1}^N|X^k_R-X^k|^p \leq 
\frac{1}{N}\sum_{k=1}^N Z^k,
$$
where $Z^k :=  2^p \, \vert X^k \vert^p \, {\bf 1}_{|X^k|>R}$ 
$(k=1, \ldots, N)$. Then, for any $p\in [1,2]$, we can introduce parameters
$\eps$ and $\theta>0$, and use Chebyshev's exponential inequality and
the independence of the variables $Z^k$ to obtain
\begin{eqnarray}
\pp\left[W_p(\hat{\mu}^N_R,\hat{\mu}^N)>\vep\right]
& \le &
\pp \left[ \frac{1}{N} \sum_{k=1}^{N} Z^k > \vep^p \right] \nonumber\\ 
& = &
\pp \left[ \exp \sum_{k=1}^{N} \theta (Z^k - \vep^p) > 1 \right] \nonumber\\
& \leq &
\ee \left( \exp \sum_{k=1}^{N} \theta(Z^k - \vep^p) \right) \nonumber\\
& = &
\exp \left( -N \left[\theta \vep^p - \ln \ee \exp(\theta Z_1) \right] \right).
\label{wp-emp-ineq}
\end{eqnarray}

In the case when \underline{$p < 2$}, for any 
$\alpha_1 < \alpha < \ds \frac{\lambda}{2}$, there exists some 
constant $R_0 = R_0(\alpha_1, p)$ such that
$$
2^p \theta r^p \leq \alpha_1 r^2 + C,
$$
for all $\theta >0$ and $r \geq R_0 \theta^{\frac{1}{2-p}}$, whence
$$
\ee \exp(\theta Z_1) \leq  \ee \exp ( \alpha_1 \, \vert X_1 \vert^2 \, {\bf 1}_{|X^k|>R} )
\leq  1 + E_{\alpha} e^{(\alpha_1 - \alpha) R^2}.
$$
As a consequence,
\begin{equation}\label{wp-emp}
\pp \left[ W_p(\hat{\mu}^N_R,\hat{\mu}^N)>\vep \right]  \, \le \, 
\exp \left(-N \left[ \theta \, \vep^p - E_{\alpha} e^{(\alpha_1 - \alpha)R^2} \right] \right).
\end{equation}

{From} \eqref{wp-trunc-ineq}, \eqref{wp-emp} and the triangular inequality for $W_p$,
\begin{eqnarray}
\pp\left[ W_p(\mu,\hat{\mu}^N)>\eps\right] 
& \le &
\pp\left[ W_p(\mu,\mu_R)+W_p(\mu_R,\hat{\mu}^N_R)+W_p(\hat{\mu}^N_R,\hat{\mu}^N)>\eps \right]\nonumber\\
&\le&
\pp\left[ W_p(\mu_R,\hat{\mu}^N_R)>\eta\, \eps -2E_{\alpha}^{1/p} R\,e^{-\frac{\alpha}{p} R^2}\right]
+ \pp\left[ W_p(\hat{\mu}^N_R,\hat{\mu}^N)>(1-\eta)\eps\right]
\nonumber\\
& \le &
\pp\left[W_p(\mu_R,\hat{\mu}^N_R)>\eta\, \eps -2E_{\alpha}^{1/p} Re^{-\frac{\alpha}{p}R^2}\right] \nonumber\\
& & + \exp \left( -  N \left(  \theta(1-\eta)^p\eps^p -
E_{\alpha}\, e^{(\alpha_1 - \alpha)\, R^2} \right) \right). \label{transi}
\end{eqnarray} 
This estimate was established for any given $p\in [1,2)$, $ \eta \in (0,1)$, $\eps, \theta >0$, 
$\ds \alpha_1 < \alpha < \frac{\lambda}{2}$ and $R \ge \max \left( \sqrt{p/{2\alpha}}, R_{0} \theta^{\frac{1}{2-p}}\right)$, 
where $R_0$ is a constant depending only on $\alpha_1$ and $p$.

\medskip

In the case when \underline{$p = 2$}, we let $Z^k :=   \vert Y_k - X_k \vert^2 \, {\bf 1}_{|X_k|>R}$ 
$(k=1, \ldots, N)$, and starting from inequality~\eqref{wp-emp-ineq} again, 
we choose $\alpha_1 < \alpha$ and then $\theta := \alpha_1/2$: 
by definition of $Z_1$ and $\mu_R$,
\begin{eqnarray*}
\ee \left( \exp\, \left( \frac{\alpha_1}{2} Z_1 \right) \right) 
& = &  
\int_{\rr^{2d}} \exp \left( \frac{\alpha_1}{2} \vert y - x \vert^2 
{\bf 1}_{\vert x \vert \geq R} \right) \, d\mu(x) \, d\mu_R(y) \\
& = & 
\mu[B_R] + \frac{1}{\mu[B_R]}   \int_{\vert y \vert \leq R} \int_{\vert x \vert \geq R} 
\exp \left( \frac{\alpha_1}{2} \vert y - x \vert^2 \right)\, d\mu(x) \, d\mu(y)  \\
& \leq &
1 + (1 - E_{\alpha} e^{-\alpha R^2})^{-1} \int_{\vert y \vert \leq R} e^{\alpha_1 \vert y \vert^2}\, d\mu(y) \, 
\int_{\vert x \vert \geq R} e^{\alpha_1 \vert x \vert^2}\, d\mu(x) \\
& \leq &
1 + 2   \, E_{\alpha}^{2} \, e^{(\alpha_1 - \alpha)R^2}
\end{eqnarray*}
for $R$ large enough, from which
\begin{equation}\label{wp-empp=2}
\pp\left[W_2(\hat{\mu}^N_R,\hat{\mu}^N)>\vep\right]  \, \le \, 
\exp \left(-N \left[ \frac{\alpha_1}{2} \, \vep^2 - 2 \, 
E_{\alpha}^2\, e^{(\alpha_1 - \alpha)R^2} \right] \right). 
\end{equation}

To sum up, in the case $p=2$ equation \eqref{transi} writes
\begin{multline}
\pp\left[ W_2(\mu,\hat{\mu}^N)>\eps\right)] \le 
\pp\left[W_2(\mu_R,\hat{\mu}^N_R)>\eta\, \eps 
-2E_{\alpha}^{1/2} Re^{-\frac{\alpha}{2}R^2}\right]
\\ + \exp \left( -  N \left(  \frac{\alpha_1}{2}(1-\eta)^2\eps^2 -
 2 \, E_{\alpha}^2  \, e^{(\alpha_1 - \alpha)\, R^2} \right) \right). \label{transip=2}
\end{multline} 

\medskip

So, apart from some error terms, for all $p\in [1,2]$ we have reduced the initial problem 
to establishing the result only for the probability law $\mu_R$, whose support lies 
in the compact set $B_R$.

\medskip

We end up this truncation procedure by proving that $\mu_R$ satisfies some modified 
$T_p$ inequality. Let indeed $\nu$ be a probability measure on 
$B_R$, absolutely continuous with respect to $\mu$ (and hence with respect to $\mu_R$); then,
when $R$ is larger than some constant depending only on $E_\alpha$, we can write
\begin{align}
H(\nu|\mu_R)-H(\nu|\mu) 
= \int_{B_R}\ln \frac{d\nu}{d\mu_R}\,d\nu - \int_{B_R}\ln \frac{d\nu}{d\mu}\, d\nu = \ln\mu[B_R] \nonumber\\
& \ge \ln\left(1-E_{\alpha} e^{-\alpha R^2}\right) \nonumber\\
& \ge -2E_{\alpha} e^{-\alpha R^2}. \label{HH}
\end{align}
But $\mu$ satisfies a $T_p(\lambda)$ inequality, so
$$
H(\nu|\mu) \ge 
\frac{\lambda}{2}\,W_p^2(\mu,\nu) \ge 
\frac{\lambda}{2}\, \Bigl(W_p(\mu_R,\nu)-W_p(\mu_R,\mu)\Bigr)^2
$$
by triangular inequality. Combining this with~\eqref{HH}, we obtain
$$
H(\nu|\mu_R) \ge \frac{\lambda}{2}\, \Bigl(W_p(\mu_R,\nu)-W_p(\mu_R,\mu)\Bigr)^2 
-2E_{\alpha}\, e^{-\alpha R^2}
$$
From this, inequality~\eqref{wp-trunc-ineq} and the elementary inequality
\begin{equation}\label{devcarre}
\forall a\in (0,1)\qquad \exists C_a>0; \qquad \forall x,y \in \rr,\qquad
(x-y)^2 \geq (1-a) \, x^2 - C_a \, y^2,
\end{equation}
we deduce that for any $\lambda_1 < \lambda$ there exists 
some constant $K$ such that
\begin{equation}\label{transp-app}
H(\nu|\mu_R)  \geq \frac{\lambda_1}{2}W_p(\mu_R,\nu)^2-KR^2e^{-\alpha R^2}.
\end{equation}

\bigskip

{\bf Step~2: Covering by small balls}. In this second step we derive quantitative 
estimates on $\hat{\mu}^N_R$. Let $\phi$ be a bounded continuous function on $\rr^d$,
and let ${\cal B}$ be a Borel set in $\PP(B_R)$ (equipped with the weak topology
of convergence against bounded continuous test functions).
By Chebyshev's exponential inequality and the independence of the variables $X^k_R$,
\begin{eqnarray*}
\ds
\pp[\hat{\mu}^N_R\in \BB]
&\le& \exp \left( -N\inf_{\nu\in\BB}\int_{B_R}\phi \, d\nu \right) 
\,\ee\left(e^{N\int_{B_R} \phi \, d\hat{\mu}^N_R}\right)\\
&=&
\exp\left(-N\inf_{\nu\in\BB}\left[\int_{B_R} \phi \, d\nu-
\frac{1}{N}\log\ee\left(e^{N\int_{B_R} \phi \, d\hat{\mu}^N_R}\right)\right]\right)\\
&=&
\exp\left(-N\inf_{\nu\in\BB}\left[\int_{B_R} \phi \, d\nu
-\frac{1}{N}\log\ee\left(e^{\sum_{k=1}^N \phi(X_R^k)}\right)\right]\right)\\
&=& \exp\left(-N\inf_{\nu\in\BB}\left[\int_{B_R} \phi \, d\nu
-\log\int_{B_R}  e^\phi d\mu_R\right]\right).
\end{eqnarray*}

As $\phi$ is arbitrary, we can pass to the supremum and find
$$
\pp[\hat{\mu}^N_R\in \BB]\le \exp\left(-N\sup_{\phi\in C_b(\rr^d)}\inf_{\nu\in\BB}
\left[\int_{B_R} \phi \, d\nu-\log\int_{B_R} e^\phi\, d\mu_R\right]\right).
$$
Now we note that the quantity $\int\phi \, d\nu-\log\int e^\phi\,d\mu_R$ is linear in 
$\nu$ and convex lower semi-continuous (with respect to the topology of
uniform convergence) in $\phi$ ; if we further assume that
$\BB$ is convex and compact, then (for instance) Sion's min-max 
theorem~\cite[Theorem 4.2']{Sio58}  ensures that
$$
\sup_{\phi\in C_b(\rr^d)} \inf_{\nu \in \BB} 
\left[\int_{B_R}\phi \, d\nu-\log\int e^\phi d\mu_R\right] = \inf_{\nu \in \BB}  \sup_{\phi\in C_b(\rr^d)}
\left[\int_{B_R} \phi \, d\nu-\log\int e^\phi d\mu_R\right].
$$
By the dual formulation of the $H$ functional~\cite[Lemma 6.2.13]{DZ98}, we conclude that
\begin{equation}\label{estBB}
\pp[\hat{\mu}^N_R\in \BB]\le \exp\left(-N\inf_{\nu\in\BB}H(\nu|\mu_R)\right).
\end{equation}

Now, let $\delta >0$ and let $\AA$ be a measurable subset of $\PP(B_R)$. We cover
the latter with $\NN^\AA$ balls $(B_i)_{1 \leq i \leq \NN^\AA}$ with radius 
$\delta/2$ in $W_p$ metric. Each of these balls is convex and compact, and it
is included in the $\delta$-thickening of $\AA$ in $W_p$ metric, defined as
$$
\AA_\delta:=\Bigl\{\nu\in\PP(B_R);\quad \exists \, \nu_a\in\AA, \; W_p(\nu,\nu_a)\le\delta\Bigr\}.
$$
So, by~\eqref{estBB} we get
\begin{eqnarray}
\pp[\hat{\mu}^N_R\in \AA]
&\le&
\pp\left[\hat{\mu}^N_R\in \bigcup_{i=1}^{\NN^\AA}B_i\right] \nonumber\\
&\le&
\sum_{i=1}^{\NN^\AA}\pp\left(\hat{\mu}^N_R\in B_i\right)\nonumber\\
&\le&
\sum_{i=1}^{\NN^\AA}\exp\left(-N\inf_{\nu\in B_i}H(\nu|\mu_R)\right) \nonumber\\
&\le&
\NN^\AA\exp\left(-N\inf_{\nu\in \AA_\delta}H(\nu|\mu_R)\right). \label{estAA}
\end{eqnarray}

We now apply this estimate with
\[ \AA:=\Bigl\{\nu\in\PP(B_R);\quad W_p(\nu,\mu_R) \ge \eta\eps - 
2 E_{\alpha}^{1/p} R \, e^{-\frac{\alpha}{p}R^2}\Bigr\}.\]
From~\eqref{transp-app} we have, for any $\nu\in\AA_\delta$,
$$
H(\nu|\mu_R) \geq \frac{\lambda_1}{2} W_p(\nu,\mu_R)^2 - K R^2 e^{-\alpha R^2} \geq
\frac{\lambda_1}{2}\rho^2 - KR^2e^{-\alpha R^2},
$$
where
$$
\rho := \max \left( \eta\, \eps - 2 E_{\alpha}^{1/p} R e^{-\frac{\alpha}{p}R^2} - \delta , 0 \right).
$$
Combining this with~\eqref{estAA}, we conclude that
\begin{equation}\label{conc-C2}
\pp\left[W_p(\mu_R,\hat{\mu}^N_R)\ge \eta\eps - 2 E_{\alpha}^{1/p} R e^{-\frac{\alpha}{p}R^2}\right]
\le
\NN^\AA \exp\left(-N\left[\frac{\lambda_1}{2}\rho^2 - KR^2e^{-\alpha R^2}\right]\right). 
\end{equation}

\medskip

Now, given any $\lambda_2 < \lambda_1$, it follows from~\eqref{devcarre} that
there exist $\delta_1$, $\eta_1$ and $K_1$, depending on $\alpha, \lambda_1, \lambda_2$, such that
\begin{equation}\label{estrho}
\frac{\lambda_1}{2}\rho^2 - KR^2e^{-\alpha R^2} \geq 
\frac{\lambda_2}{2}\eps^2 - K_1 R^2e^{-\alpha R^2}
\end{equation}
where $\delta := \delta_1 \eps$ and $\eta := \eta_1$.

Though this inequality holds independently of $p$, we shall use it only in the case when $p < 2$.
In the case $p=2$, on the other hand, we note that for any $\eta\in (0,1)$,
\begin{equation}\label{estrhop=2}
\frac{\lambda_1}{2}\rho^2 - KR^2e^{-\alpha R^2} \geq \frac{\lambda_2}{2} \eta^2\eps^2 
- K_1 R^2e^{-\alpha R^2}
\end{equation}
where $\delta := \delta_1 \eps$.

Finally, we bound $\NN^\AA$ by means of Theorem~\ref{Th-balls} in 
Appendix \ref{sectionentropiemetrique}: there exists some constant $C$ (only depending on $d$)
such that for all $R >0$ and $\delta > 0$ the set $\PP(B_R)$ can be covered by
\[
\left( C\frac{R}{\delta} \vee 1\right)^{\left(C\frac{R}{\delta}\right)^d} 
\]
balls of radius $\delta$ in $W_p$ metric, where $a \vee b$ stands for $\max(a,b)$.
In particular, given $\delta = \delta_1 \eps$, we can choose
\begin{equation}\label{estNNAA}
\NN^\AA \leq  \left( K_2\frac{R}{\eps} \vee 1 \right)^{\left(K_2\frac{R}{\eps}\right)^d}
\end{equation}
balls of radius $\delta$,  for some constant $K_2$ depending on $\lambda_1$ and $\lambda_2$ 
(via $\delta_1$) but neither on $\eps$ nor on $R$.
(The purpose of the~1 in $(K_2R/\var\vee 1)$ is to make sure that
the estimate is also valid when $\var>R$.)
\med

Combining~\eqref{transi}, \eqref{conc-C2}, \eqref{estrho} and \eqref{estNNAA}, we find that,
given {\bf $p \in [1,2)$}, $\lambda_2 < \lambda$ and $\ds \alpha_1 < \alpha < \frac{\lambda}{2}$, 
there exist some constants $K_1$, $K_2$, $K_3$ and $R_1$ such that for all 
$\eps , \zeta >0$ and $R \geq R_1 \max (1, \zeta^{\frac{1}{2-p}})$,
\begin{multline} \label{horreurp<2}
\pp\left[W_p(\mu,\hat{\mu}^N) > \eps \right]
\le \left( K_2\frac{R}{\eps} \vee 1 \right)^{K_2\left(\frac{R}{\eps}\right)^d}
\exp\left(-N\left[\frac{\lambda_2 \eps^2}{2} - K_1R^2e^{-\alpha R^2}\right]\right) \\
+ \exp \left( -  N \left(K_3 \, \zeta \, \eps^p - 
K_4 e^{(\alpha_1 - \alpha)\, R^2} \right) \right) 
\end{multline}
for some constant $K_4 = K_4(\theta, \alpha_1)$. In the case when 
\underline{$p=2$}, we obtain similarly
\begin{multline}
\pp\left[W_2(\mu,\hat{\mu}^N) > \eps \right]
\le 
\left( K_2\frac{R}{\eps} \vee 1 \right)^{K_2\left(\frac{R}{\eps}\right)^d}
\exp\left(-N\left[\frac{\lambda_2}{2} \eta^2 \eps^2 - 
K_1R^2e^{-\alpha R^2}\right]\right) \\
+ \exp \left( -  N \left(\frac{\alpha_1}{2} (1 - \eta)^2  \, \eps^2 - 
K_4 e^{(\alpha_1 - \alpha)\, R^2} \right) \right)\label{conc-C2'} 
\end{multline}
for any $\eta \in (0,1)$ and $R \geq R_1$.

These estimates are not really appealing (!), but they are rather precise and general. In the rest
of the section we shall show that an adequate choice of $R$ leads to a simplified expression.
\med

{\bf Step~3: Choice of the parameters}.

We first consider the case when \underline{$p \in [1,2)$}. Let $\lambda'<\lambda_2$, $\alpha' < \alpha$ 
and $d_1 >d$. We claim that
$$
\pp\left[W_p(\mu,\hat{\mu}^N) > \eps \right]
\le \exp \left(-\frac{\lambda'}{2} N \, \eps^2\right) + \exp \, (- \alpha' N \, \eps^2)
$$
as soon as
\begin{equation}\label{conditions}
           R^2  \geq R_2 \max \left( 1 , \eps^2, \ln \left(\frac{1}{\eps^2} \right) \right) , \qquad
           N \, \eps^{d_1+2} \geq K_5 R^{d_1} 
\end{equation}
for some constants $R_2$ and $K_5$ depending on $\mu$ only through $\lambda, \alpha$ 
and $E_{\alpha}$.

Indeed, on one hand
$$
K_2 \left( \frac{R}{\eps}\right)^d \ln \left( K_2 \frac{R}{\eps}\right) 
\leq K_6 \left(\frac{R}{\eps}\right)^{d_1}
$$
for some constant $K_6$, on the other hand
$$
K_1R^2e^{-\alpha R^2} \leq e^{-\alpha_1 R^2}
$$
for $R$ large enough, and then
$$
K_6 \left(\frac{R}{\eps}\right)^{d_1} - N \left[\frac{\lambda_2 \eps^2}{2} - 
e^{-\alpha_1 R^2} \right] \leq -N \, \frac{\lambda' \eps^2}{2}
$$
for $R^2 / \ln (\frac{1}{\eps^2})$ and $N \, \eps^{d_1+2} / R^{d_1}$ large enough; this
is enough to bound the first term in the right-hand side of~\eqref{horreurp<2} if moreover
$R/\eps$ is large enough. 

Moreover, letting $\alpha_2 \in (\alpha', \alpha_1)$, we can choose $\zeta$ in such 
a way that $ K_3 \zeta = \eps^{2-p}$, so that
$$
\exp \left( - N \left( K_3 \zeta \eps^p - K_4 \, e^{(\alpha_1 - \alpha)\, R^2} \right) \right)
= \exp \left(- N \left( \alpha_2 \, \eps^2 - K_4 e^{(\alpha_1 - \alpha)\, R^2} \right) \right),
$$
which in the end can be bounded by
$$
\exp \, ( - N \, \alpha' \eps^2)
$$
if $R$ and $R^2 / \ln(\frac{1}{\eps^2})$ are large enough. With this one can get a bound on
the right-hand side of~\eqref{horreurp<2}.
\medskip

Now let us check that conditions~\eqref{conditions} can indeed be fulfilled. 
Clearly, the first condition holds true for all $\eps \in (0,1)$ and 
$R^2 \geq R_3 \ln (\frac{K_6}{\eps^2})$, where $R_3$ and $K_6$ are positive constants. 
Then, we can choose
$$
R := \left( \frac{N}{K_5} \eps^{d_1+2} \right)^{1/d_1}
$$
so that the second condition holds as an equality. This choice is admissible as soon as
$$
\left( \frac{N}{K_5} \, \eps^{d_1+2} \right)^{2/d_1} \geq R_3 \ln \left(\frac{K_5}{\eps^2}\right)
$$
and this, in turn, holds true as soon as
\begin{equation}\label{condneps}
N \geq K_7 \, \, \eps^{-(d'+2)},
\end{equation}
where $d'$ is such that $d'>d$, and $K_7$ is large enough.
\med

If $\var \geq1$, then we can choose $R^2=R_2 \var^2$, i.e. $R=\sqrt{R_2} \var$, and then
the second inequality in~\eqref{conditions} will be true as soon as $N$ is large enough.
\med

To sum up: Given $d'>d$, $\lambda' < \lambda$ and $\alpha' < \alpha$, 
there exists some constant $N_0$, depending on $d'$ and depending on $\mu$ 
only through $\lambda, \alpha$ and $E_{\alpha}$, such that for all $\eps>0$,
$$
\pp\left[W_p(\mu,\hat{\mu}^N) > \eps \right] \le  
\exp \left(- \frac{\lambda'}{2} \, N \, \eps^2\right) 
+ \exp \left(- \alpha' \, N \, \eps^2\right)
$$
as soon as $N \geq N_0 \max(\eps^{-(d'+2)},1)$. 
Then we note that, given $\ds K <  \min \left( \frac{\lambda'}{2}\, , \alpha' \right)$, 
the inequality
$$
\exp \left(- \frac{\lambda'}{2} \, N \, \eps^2\right) + \exp \left(-  \alpha' \, N \, \eps^2\right)
\leq \exp \left(- K \, N\, \eps^2\right)
$$
holds if condition~\eqref{condneps} is satisfied for some $K_7$ large enough. To conclude
the proof of Theorem~\ref{thmconc-Tp2} in the case when $p\in [1,2)$, it is sufficient
to choose $\lambda'<\lambda$, $\alpha<\lambda/2$.

\bigskip
Now, in the case when \underline{$p = 2$}, given $\lambda_3 < \lambda_2$ and $\alpha_2 < \alpha_1$, 
conditions~\eqref{conditions} imply
$$
\pp\left[W_2(\mu,\hat{\mu}^N) > \eps \right]
\le \exp \left(-\frac{\lambda_3}{2} \eta^2 \, N\eps^2\right) + 
\exp \, \left(-\frac{\alpha_2}{2} (1-\eta)^2\, N\eps^2\right).
$$
Then we let $\ds \alpha_2 := \frac{\lambda_3}{2}$ and $\eta := \sqrt{2} - 1$, so that
$$
\frac{\lambda_3}{2} \, \eta^2 = \frac{\alpha_2}{2} \, (1-\eta)^2.
$$
Then
$$
\pp\left[W_2(\mu,\hat{\mu}^N) > \eps \right] \leq 2 \, 
\exp \left(-(3 - 2 \, \sqrt{2}) \frac{\lambda_3}{2} N  \, \eps^2\right);
$$
for $\lambda' < \lambda$, the above quantity is bounded by
$$
\exp \left(-(3 - 2 \, \sqrt{2}) \frac{\lambda'}{2} N  \, \eps^2\right)
$$
as soon as~\eqref{condneps} is enforced with $K_7$ large enough. 
This concludes the argument.

\subsection{Proof of Theorem~\ref{thmconc-Mq}}

It is very similar to the proof of Theorem~\ref{thmconc-Tp2}, so we shall only explain
where the differences lie. Obviously, the main difficulty will consist in the control of tails.

We first let $p \in [1,q)$, $\ds \alpha \in [1, \frac{q}{p})$ and $R >0$, 
and introduce 
\[ M_q := \int_{\rr^d} \vert x \vert^q \, d\mu(x).\] 
Then~\eqref{wp-trunc-ineq} may be replaced by
\begin{equation}\label{wp-trunc-ineq-pol}
W_p^p(\mu,\mu_R) \leq 2^p M_q R^{p-q},
\end{equation}
and~\eqref{wp-emp-ineq} by
\begin{equation}\label{wp-emp-ineq-pol}
\pp \, [W_p(\hat{\mu}^N_R,\hat{\mu}^N) > \vep] \leq C \, N^{{\overline{\alpha}} - \alpha} 
\frac{R^{\alpha \, p -q}}{(\vep^p - C \, R^{p-q})^{\alpha}}
\end{equation}
for some constant $C$ depending on $\alpha$ and $M_q$.

Let us establish for instance~\eqref{wp-emp-ineq-pol}. Introduce
\[ Z_k = \vert Y_k - X_k \vert^{p} \, {\bf 1}_{\vert X_k \vert > R} \qquad
(1\leq k\leq N). \]
By Chebychev's inequality,
\begin{multline*}
\pp \, \left[W_p(\hat{\mu}^N_R,\hat{\mu}^N) > \vep\right] \leq 
\pp \, \left[ \frac {1}{N} \sum_{k=1}^{N} Z_k > \vep^p \right]
= \pp \, \left[ \frac{1}{N} \sum_{k=1}^{N} (Z_k - \ee \, Z_k) > \vep^p - \ee \, Z_1  \right] \\
\leq \frac{ \ee \left\vert \sum_{k=1}^{N} (Z_k - \ee \, Z_k) \right\vert^{\alpha}}
{(N \, (\vep^p - \ee \, Z_1))^{\alpha}}
\end{multline*}
provided that $\vep^p > \ee \, Z_1$. But, since the random variables 
$(Z_k - \ee \, Z_k)_k$ are independent and identically distributed, with zero mean, 
there exists some constant $C$ depending on $\alpha$ such that
$$
\ee \left\vert \sum_{k=1}^{N} ( Z_k - \ee \, Z_k ) \right\vert^{\alpha} \leq 
C \, N^{\overline{\alpha}} \, \ee \vert Z_1 - \ee \, Z_1 \vert^{\alpha}
$$ 
where $\ds {\overline{\alpha}} := \max ( \alpha/2 , 1)$. 
This inequality is a consequence of Rosenthal's inequality in the case 
when $\alpha \geq 2$, but also holds true if $\alpha \in [1,2)$ 
(see for instance~\cite[pp.~62 and~82]{P95}). Then, on one hand,
$$
\ee \, Z_1 = \ee \, \vert Y_1 -  X_1 \vert^{p}  \, {\bf 1}_{\vert X_1 \vert >R} \leq 2^p M_q R^{p-q},
$$
while on the other hand,
\begin{multline*}
\ee \vert Z_1 - \ee Z_1 \vert^{\alpha} 
= \ee \left\vert \vert Y_1 -  X_1 \vert^{p}  \, {\bf 1}_{\vert X_1 \vert >R}
- \ee \vert Y_1 -  X_1 \vert^{p}  \, {\bf 1}_{\vert X_1 \vert >R} \right\vert^{\alpha}
\\ \leq C \ee \vert Y_1 -  X_1 \vert^{\alpha p}  \, {\bf 1}_{\vert X_1 \vert >R}
\leq C \, M_q \, R^{\alpha p -q}
\end{multline*}
with $C$ standing for various constants. Collecting these two estimates, we conclude to
the validity of~\eqref{wp-emp-ineq-pol} for $R^{q-p} \, \vep^p$ large enough.

\medskip

Then~\eqref{wp-trunc-ineq-pol} and~\eqref{wp-emp-ineq-pol} together ensure that
\begin{multline}\label{transiMq}
\pp[W_p(\mu,\hat{\mu}^N)> \eps] \leq 
\pp\left[ W_p(\mu_R,\hat{\mu}^N_R)>\eta\, \eps -2M_q^{1/p} R^{1-q/p} \right]
+ \, C \, N^{{\overline{\alpha}} - \alpha} \frac{R^{\alpha \, p -q}}
{((1 - \eta)^p \eps^p - C \, R^{p-q})^{\alpha}}
\end{multline}
for any $\eps \in (0,1)$, $\eta > 0$ and $R^{q-p} \, \eps^p \, (1-\eta)^p$ large enough.

Since $\mu_R$ is supported in $B_R$, the Csisz\'ar-Kullback-Pinsker inequality and 
Kantorovich-Rubinstein formulation of the $W_1$ distance together ensure that it 
satisfies a $T_1(R^{-2})$ inequality (see e.g.~\cite[Particular Case~5]{BV04} with $p=1$).  
This estimate also extends to any $W_p$ distance, not as a penalized
$T_p$ inequality as in~\eqref{transp-app}, but rather as
\begin{equation}\label{modTp}
W_p^{2p}(\nu,\mu_R) \leq 2^{2p-1} R^{2p} H(\nu \vert \mu_R)
\end{equation}
(see again~\cite[Particular Case 5]{BV04}).

From~\eqref{transiMq} and~\eqref{modTp} we deduce (as in~\eqref{conc-C2'}) that
\begin{multline}\label{conc-Cq}
\pp\left[W_p(\hat{\mu}^N,\mu) > \eps \right]
\le
\left( K_1\frac{R}{\delta} \right)^{K_1\left(\frac{R}{\delta}\right)^d}
 \exp\left(-\frac{N \rho^{2p}}{2^{2p-1}R^{2p}} \right) 
+ C \, N^{{\overline{\alpha}} - \alpha} \frac{R^{\alpha \, p -q}}
{((1 - \eta)^p\eps^p - C \, R^{p-q})^{\alpha}}
\end{multline}
for any $\delta$, where now
$$
\rho := \left( \eta \eps - 2M^{1/p}R^{1-q/p} - \delta \right)^+.
$$

Letting $\eta_1 < \eta$ and $d' > d$, and choosing $\delta = \delta_0 \, \eps$, we deduce
$$
\pp\left[(W_p(\hat{\mu}^N,\mu) > \eps \right]
\le
\exp\left( \left(\frac{R}{\eps}\right)^{d'} -
\frac{\eta_1^{2p}}{2^{2p-1}} \frac{N \, \eps^{2p}}{R^{2p}} 
+ \frac{K_1}{2^{2p-1}} \frac{N}{R^{2q}} \right) 
+ C \, N^{\bar\alpha - \alpha} \frac{R^{\alpha \, p -q}}
{((1 - \eta_1)^p\eps^p - C \, R^{p-q})^{\alpha}}
$$
for $R^{q-p} \, \eps^p \, (1-\eta_1)^p$  large enough, and then
\begin{equation}\label{conc-Cq'}
\pp\left[W_p(\hat{\mu}^N,\mu) > \eps \right]
\le
\exp\left( -\frac{\eta_2^{2p}}{2^{2p-1}} \frac{N \, \eps^{2p}}{R^{2p}}  \right) 
+ C \, N^{{\overline{\alpha}} - \alpha} \frac{R^{\alpha \, p -q}}
{(1 - \eta_2)^{\alpha p}\eps^{\alpha p}}
\end{equation}
for $\eta_2 < \eta_1$, provided that the conditions
\begin{equation}\label{conditions'}
           R  \geq R_1  \eps^{-\frac{p}{q-p}}, \qquad
           N  \geq K_2 \left( \frac{R}{\eps} \right)^{2p+d'}  
\end{equation}
hold for some $R_1$ and $K_2$. 

Given any choice of $R$ as a product of powers of $N$ and $\eps$, 
the first term in the right-hand side of~\eqref{conc-Cq'} will always be
smaller than the second one, if $N$ goes to infinity while $\eps$ is kept fixed; thus
we can choose $R$ minimizing the second term  under the above conditions. 
Then the second condition in~\eqref{conditions'} will be fulfilled as an equality:
$$
R = K_3 \, \eps \, N^{\frac{1}{2p+d'}}.
$$

As for the first condition in~\eqref{conditions'}, it can be rewritten as
$$
N \geq N_0 \, \eps^{-q \frac{2p+d'}{q-p}},
$$
and then, by \eqref{conc-Cq'},
$$
\pp \left[W_p(\hat{\mu}^N,\mu) > \eps \right]
 \le 
\exp \left( - K_5 N^{\frac{d'}{2p+d'}} \right) 
+ K_6 \, \eps^{-q} \, N^{{\overline \alpha} - \alpha + \frac{\alpha p-q}{2p+d'}}.
$$
Hence
\begin{equation}\label{alpha}
\pp \left[W_p(\hat{\mu}^N,\mu) > \eps \right] \le \eps^{-q}  \, N^{{\overline{\alpha}} - \alpha}
\end{equation}
for all $\eps \in (0,1)$ and $N$ larger than some constant and, given $d' >d$, for all
$\eps \geq 1$ and $N \geq M \eps^{d' -d}$ where $M$ is large enough.

\medskip

In the first case when $p  \geq q/2$, any admissible $\alpha$ belongs to 
$[1,q/p) \subset [1,2]$, so $\overline{\alpha} = 1$. If $\delta \in (0, q/p - 1)$, 
we get from~\eqref{alpha}, with $\alpha = q/p - \delta$, that
$$
\pp\Bigl[W_p(\hat{\mu}^N,\mu) > \eps \Bigr]
\le
\eps^{-q} N^{1 - q/p + \delta} 
$$
for all $\eps > 0$ and 
$$
N \geq N_0 \max \bigl(\eps^{-q \frac{2p+d'}{q-p}} , \eps^{d' -d} \bigr).
$$

In the second case when $p < q/2$, we only consider admissible $\alpha$'s 
in $[2,q/p) \subset [1, q/p)$, so that $\overline{\alpha} - \alpha = - \alpha /2$. 
Choosing $\delta \in (0, q/p - 2)$, we get from~\eqref{alpha}
$$
\pp\Bigl[W_p(\hat{\mu}^N,\mu) > \eps \Bigr]
\le \eps^{-q} N^{-q/2p + \delta/2} 
$$
under the same conditions on $N$ as before. This concludes the argument.

\subsection{Proof of Theorem~\ref{thmconc-var}}

It is again based on the same principles as the proofs of Theorems~\ref{thmconc-Tp2}
and~\ref{thmconc-Mq}, with the help of functional inequalities investigated
in~\cite{BV04} and~\cite{CG04}. We skip the argument, which the reader can easily
reconstruct by following the same lines as above.

\subsection{Data reconstruction estimates}

Finally, we show how the above concentration estimates imply data reconstruction estimates.
This is a rather general estimate, which is treated here along the lines
of ~\cite[Section 5]{Sch96} and~\cite[Problem~10]{V03}.

\begin{proposition} \label{propreconstr}
Let $\mu$ be a probability measure on $\rr^d$, with density $f$ with
respect to Lebesgue measure. Let $X_1,\ldots, X_N$ be random points in $\rr^d$,
and let $\zeta$ be a Lipschitz, nonnegative kernel with unit integral. 
Define the random measure $\hat{\mu}$ and the random function $\hat{f}_{\zeta,\alpha}$ by
\[ \hat{\mu} := \frac1N \sum_{i=1}^N \delta_{X_i}, \qquad
\hat{f}_{\zeta,\alpha} (x) := \frac1N \sum_{i=1}^N 
\zeta_\alpha(x-X_i), \qquad
\zeta_\alpha(x) = \frac1{\alpha^d} \zeta\left( \frac{x}{\alpha} \right).\]
Then, 
\begin{equation} \label{superror}
\sup_{x\in\R^d} |\hat{f}_{\zeta,\alpha}(x) - f(x) | \leq \frac{\|\zeta\|_\Lip}
{\alpha^{d+1}} W_1(\hat{\mu},\mu) +\delta(\alpha),
\end{equation}
where $\delta$ stands for the modulus of continuity of $f$, defined as
\[ \delta(\var) := \sup_{|x-y|\leq \var} |f(x)-f(y)|. \]
\sm

As a consequence, if $f$ is Lipschitz, then there exist some constants $a, K>0$, only depending 
on $d$, $\|f\|_\Lip$ and $\|\zeta\|_\Lip$, such that 
\begin{equation}\label{Liperror} 
\pp \Bigl[ \|\hat{f}_{\zeta,a\var} - f \|_{L^\infty} > \var\Bigr]
\leq \pp \Bigl[ W_1(\hat{\mu},\mu)> K \var^{d+2} \Bigr] 
\end{equation}
for all $\eps  >0$.
\end{proposition} 

\begin{proof} First,
\begin{align*} |\mu\ast \zeta_\alpha(x) - f (x)| =
& \left | \int_{\rr^d} \zeta_\alpha(x-y) \,\bigl(f (y)- f(x)\bigr) \,dy \right | \\
& \leq \int_{\rr^d} \zeta_\alpha(x-y) |f(y)-f(x)|\,dy.
\end{align*}
Since $\zeta_\alpha(x-y)$ is supported in $\{|x-y|\leq \alpha\}$, and
$\zeta_\alpha$ is a probability density, we deduce
\begin{equation}\label{mumuast} 
|\mu\ast \zeta_\alpha(x) - f (x)| \leq \delta(\alpha).
\end{equation}

Now, if $x$ is some point in $\rr^d$, then, thanks to the
Kantorovich-Rubinstein dual formulation~\eqref{KR},
\begin{align*} 
\Bigl |\hat{f}_{\zeta,\alpha} - \mu\ast \zeta_\alpha \Bigr |(x)
& = \left | \int_{\rr^d} \zeta_\alpha(x-y)\, d[\hat{\mu}-\mu](y) \right |\\
& \leq \|\zeta_\alpha(x-\cdot)\|_\Lip W_1(\hat{\mu}, \mu) \\
& = \frac{\|\zeta\|_\Lip}{\alpha^{d+1}} W_1(\hat{\mu}, \mu).
\end{align*}
To conclude the proof of~\eqref{superror}, it suffices to combine this bound 
with~\eqref{mumuast}.
\med

Now, let $L:=\max (\|f\|_\Lip, \|\zeta\|_\Lip)$, and $\alpha:=\var/(2L)$.
The bound \eqref{superror} turns into
\[ \|\hat{f}_{\zeta,\alpha} - f \|_{L^\infty} \leq L \bigl(\frac{W_1(\hat{\mu},\mu)}{\alpha^{d+1}} + \alpha\bigr)
\leq \left (\frac{(2L)^{d+1}L}{\var^{d+1}}\right ) W_1(\hat{\mu},\mu) + \frac{\var}2. \]
In particular, 
\[ \pp \Bigl[\|\hat{f}_{\zeta,\alpha} -f\|_{L^\infty}>\var \Bigr]
\leq\pp \left[ W_1(\hat{\mu},\mu) > \frac{\var^{d+2}}{(2L)^{d+2}} \right],\]
which is estimate~\eqref{Liperror}.
\end{proof}

\begin{remark}
Estimate \eqref{Liperror}, combined with Theorem \ref{thmconc-Tp2} or Theorem \ref{thmconc-Mq}, yields
simple quantitative (non-asymptotic) deviation inequalities for empirical distribution functions in supremum norm.
We refer to Gao \cite{Gao03} for a recent study of deviation inequalities for empirical distribution functions, both
in moderate and large deviations regimes.
\end{remark}

\section{PDE estimates}

Now we start the study of our model system for interacting particles.
The first step towards our proof of Theorem~\ref{thmconc-eds} consists in deriving
suitable a priori estimates on the solution to the nonlinear limit partial differential
equation~\eqref{edp}. In this section, we recall some estimates which have already
been established by various authors, and derive some new ones. All estimates will be effective.

\subsection{Notation} \label{subnotation}

In the sequel, $\mu_0$ is a probability measure, taken as an initial
datum for equation~\eqref{edp}, and various regularity assumptions will later be
made on $\mu_0$. Assumptions~\eqref{D2V,D2W,V} will always be made on $V$ and $W$,
even if they are not recalled explicitly; we shall only mention additional
regularity assumptions, when used in our estimates. Moreover, we shall write
\begin{equation}\label{Gamma} \Gamma := \max (\vert \gamma \vert , \vert \gamma' \vert).
\end{equation}
The notation $\mu_t$ will always stand for the solution (unique under our assumptions)
of~\eqref{edp}.

We also write
\[ e(t) := \int_{\rr^d} |x|^2\,d\mu_t(x)\]
for the (kinetic) energy associated with $\mu_t$, and
\[ M_\alpha(t) := \int_{\rr^d} e^{\alpha \vert x \vert^2} \, d\mu_t(x)\]
for the square exponential moment of order $\alpha$. 

The scalar product between two vectors $v,w\in\R^d$ will be denoted by $v\cdot w$.
The symbols $C$ and $K$ will often be used to denote various positive constants;
in general what will matter is an upper bound on constants denoted $C$, and a lower
bound on constants denoted $K$. The space ${\cal C}^k$ is the space of $k$ times
differentiable continuous functions.

\subsection{Decay at infinity}

In this subsection, we prove the propagation of strong decay estimates at infinity:

\begin{proposition}\label{propT_1edp}
With the conventions of Subsection~\ref{subnotation},
let $\overline{\eta}$ be $ - \gamma$ if $\gamma<0$, and an arbitrary negative number otherwise.
Let 
\[ a:= 2( \beta + \overline{\eta}), \qquad
\overline{G} := 2\, d + \frac{|\nabla V(0)|^2}{ 2\, |\overline{\eta}|}. \] 
Then
\sm

(i) $\ds e(t) \leq e^{-a t} \left[ e(0) + \overline{G} \, \frac{e^{a t}-1}{a} \right]$;
\sm

(ii) For any $\alpha_0>0$ there is
a continuous positive function $\alpha(t)$ such that $\alpha(0)=\alpha_0$ and
\begin{equation}\label{propalpha}
M_{\alpha_0}(0) < +\infty \Longrightarrow M_{\alpha(t)} (t) <+\infty.
\end{equation}
\sm

(iii) Moreover, in the ``uniformly convex case'' when $\beta > 0$ 
and $\beta+\gamma>0$, then there 
is $\alpha >0$ such that
\[ \sup_{t\geq 0}\: e(t) <+\infty, \qquad \sup_{t\geq 0}\: M_\alpha(t) <+\infty.\]
\end{proposition}

\begin{corollary} \label{sqexpT1} If $\mu_0$ admits a finite square exponential moment, then
$\mu_t$ satisfies $T_1(\lambda_t)$, for some function $\lambda_t>0$, bounded below
on any interval $[0,T]$ ($T<\infty$).
\end{corollary}

\begin{proof} We start with (i). For simplicity we shall pretend that $\mu_t$ is
a smoothly differentiable function of $t$, with rapid decay, so that all computations
based on integrating equation~\eqref{edp} against $|x|^2$ are justified. These assumptions
are not a priori satisfied, but the resulting bounds can easily be rigorously justified
with standard but tedious approximation arguments. With that in mind, we compute
$$
e'(t) = 2\, d - 2 \int_{\rr^d} ( x \, \cdot \, \nabla V(x) + x\, \cdot \, \nabla W \ast \mu_t(x) ) \, 
d\mu_t(x)
$$
with
$$
- 2 \int_{\rr^d} x \, \cdot \, \nabla V(x) \, d\mu_t(x) \leq  -2 \beta \int_{\rr^d} \vert x \vert^2 \, d\mu_t(x)
- 2 \nabla V(0) \, \cdot \, \int_{\rr^d} x \, d\mu_t(x).
$$
Since $\nabla W$ is an odd function, we have
\begin{eqnarray*}
-2 \int_{\rr^d} x \, \cdot \, \nabla W \ast \mu_t(x) \, d\mu_t(x)
& = & 
-2 \iint x \, \cdot \, \nabla W(x-y) \, d\mu_t(y) \, d\mu_t(x) \\
& = & 
- \iint (x-y) \, \cdot \, \nabla W(x-y) \, d\mu_t(y) \, d\mu_t(x) \\ 
& \leq &
- \gamma \iint \vert x -y \vert^2 \, d\mu_t(y) \, d\mu_t(x) \\ 
& = &
- 2 \, \gamma \left[ \int \vert x \vert^2 \, d\mu_t(x) \, 
- \, \left\vert \int  x  \, d\mu_t(x) \right\vert^2 \right].
\end{eqnarray*}

If $\gamma < 0$, then
\begin{eqnarray*}
e'(t) 
& \leq &
2\, d - 2(\gamma + \beta) e (t) + 2 \, \gamma \left\vert \int  x  \, d\mu_t(x) 
+ \frac{\nabla V(0)}{2 \, |\gamma|} \right\vert^2
+ \frac{\vert \nabla V(0) \vert^2}{2\, |\gamma|} \\
& \leq &
2\, d - 2(\gamma + \beta) e (t) + \frac{\vert \nabla V(0) \vert^2}{2\, |\gamma|},
\end{eqnarray*}

and if $\gamma \geq 0$, then for any ${\overline \eta} <0$
\begin{eqnarray*}
e'(t) 
& \leq &
2\, d - 2({\overline \eta} + \beta) e (t) - 2 \gamma \left( \int  \vert x \vert^2 \, d\mu_t(x) - 
\left\vert \int  x  \, d\mu_t(x)  \right\vert^2 \right)
+ \frac{\vert \nabla V(0) \vert^2}{2\, |{\overline \eta}|} \\
& \leq &
2\, d - 2( {\overline \eta} + \beta) e (t) + 
\frac{\vert \nabla V(0) \vert^2}{2\, |{\overline \eta}|} \cdot
\end{eqnarray*}

This leads to
$$
e'(t) \leq \overline{G} - a \, e(t),
$$
and the conclusion follows easily by Gronwall's lemma.

\med

We now turn to (ii). Let $\alpha$ be some arbitrary nonnegative 
${\mathcal C}^1$ function on $\rr_+$. By using the equation~\eqref{edp}, we compute
$$
\frac{d}{dt} \int e^{\alpha(t) \vert x \vert^2} \, d\mu_t(x) = 
\int \bigl[ 2 d \alpha + 4 \alpha^2 \vert x \vert^2 - 2 \alpha x \cdot \nabla V(x) 
- 2 \alpha x \cdot \nabla W\ast\mu_t(x) 
+ \alpha'(t) \vert x \vert^2 \bigr] e^{\alpha \vert x \vert^2} \, d\mu_t(x).
$$
Since $D^2V(x) \geq \beta I$ for all $x \in \rr^d$, we can write
\begin{equation}\label{Vexp}
-x \, \cdot \, \nabla V(x) \leq  -x \, \cdot \, \nabla V(0) - \beta \vert x \vert^2 \leq
- \beta \vert x \vert^2 + |\nabla V(0)| \vert x \vert
\leq (\delta - \beta) \vert x \vert^2 + 
\frac{C}{4 \, \delta}
\end{equation}
for any $\delta > 0$ and $x \in \rr^d$.

Next, our assumptions on $W$ imply $\nabla W(0)=0$, and $ \gamma I \leq D^2W(x) \leq \gamma' I$,
so
$$
x \, \cdot \, \nabla W(x) \geq \gamma \, \vert x \vert^2 \qquad \text{and}
 \qquad \vert x \, \cdot \, D^2 W(z) \, y \vert \leq \Gamma \vert x \vert
\, \vert y \vert
$$
for all $x, y, z \in \rr^d$, with $\Gamma$ defined by~\eqref{Gamma}.
Hence, by Taylor's formula,
\begin{eqnarray}
-x \, \cdot \, \nabla W\ast \mu_t(x)
& = &
- \int_{\rr^d} x \, \cdot \, \nabla W(x-y) \, d\mu_t(y) \nonumber \\
& = &
- x \, \cdot \, \nabla W(x) + \int_{\rr^d} \int_{0}^{1} x \, \cdot \, 
D^2W(x-s\, y)\, y \, d\mu_t(y) \, ds \nonumber \\
& \leq & 
-\gamma \, \vert x \vert^2 + 
\Gamma\, \vert x \vert \int_{\rr^d} \vert y \vert \, d\mu_t(y) \nonumber \\
& \leq & (-\gamma + \Gamma \eta) \vert x \vert^2 \, + \,  \frac{\Gamma}{4 \, \eta}\, 
e(t), \label{Wexp}
\end{eqnarray}
where $\eta$ is any positive number.

From \eqref{Vexp} and \eqref{Wexp} we obtain
\begin{equation}\label{derexp}
\frac{d}{dt} \Bigl( M_{\alpha(t)}(t) \Bigr) \leq
\int_{\rr^d} [A(t) + B(t) \vert x \vert^2 ] \, e^{\alpha(t) \vert x \vert^2} \, d\mu_t(x)
\end{equation}
where
\[ A(t) = C \alpha(t)\, (1 + e(t)), \qquad B(t) = \alpha'(t) + 4 \, \alpha(t)^2 + b \, \alpha(t),\]
and $C$ is a finite constant, while $b = -2 (\gamma + \beta - \delta -\Gamma \eta)$.

We now choose $\alpha(t)$ in such a way that $B(t)\equiv 0$, i.e.
$$
\alpha'(t) + 4 \, \alpha ^2(t) + b \, \alpha(t) = 0, \qquad \alpha(0)=\alpha_0.
$$
This integrates to
$$
\alpha(t) = e^{-b\, t} \left( \frac{1}{\alpha_0} + 4 \, \frac{1 - e^{-b \, t}}{b} \right)^{-1}
\quad \left( = \, \left( \frac{1}{\alpha_0} + 4 \, t \right)^{-1} \; \text{ if } \, b=0 \right).
$$
Obviously $\alpha$ is a continuous positive function, and our estimates imply
$$
\frac{d}{dt} \Bigl( M_{\alpha(t)}(t) \Bigr) \leq A(t) M_{\alpha(t)}(t).
$$
We conclude by using Gronwall's lemma that
$$
M_{\alpha(t)} (t) \leq \exp \left( \int_0^t A(s) \, ds \right) \, M_{\alpha_0}(0).
$$

\med

Next, the estimate (iii) for $e(t)$ is an easy consequence of our explicit estimates when 
$\beta >0, \beta + \gamma > 0$ (in the case when $\gamma\geq 0$ and $\beta>0$, we choose 
$\overline{\eta}\in (0,\beta)$).
%; when the center of mass is kept invariant by the equation, 
%estimate (i) holds with
%\[
%a:= 2( \beta + \gamma), \qquad
%\overline{G} := 2\, d - 2 \, \nabla V(0) \int x \, d\mu_t(x) + 2 \, \gamma 
%\left\vert \int x \, d\mu_t(x) \right\vert^2,
%\] 
%so that estimate (iii) follows if moreover $\beta + \gamma >0$.

As for the estimate about $M_\alpha(t)$, it will result from a slightly more precise
computation. 
From~\eqref{derexp}, we have
\begin{equation}\label{derexpconvexe}
\frac{d}{dt} \int_{\rr^d} e^{\alpha \, \vert x \vert^2} \, d\mu_t(x) \leq 
 \int_{\rr^d} [A(t) + B \vert x \vert^2] \, e^{\alpha \vert x \vert^2} \, d\mu_t(x)
\end{equation}
where $A$ is bounded on $\rr_+$ by some constant $a$, and
$$
B = 2 \, \alpha \, [ 2 \alpha - (\beta+\gamma - \delta -\Gamma \eta)].
$$
Since $\beta  + \gamma >0$, for any fixed $\alpha $ in 
$\ds \left( 0 \, , \frac{\beta +\gamma}{2} \right)$ we can 
choose $\delta, \eta > 0$ such that $B < 0$. Letting $R^2 = -a / B$ and $G = -B >0$, 
equation~\eqref{derexpconvexe} becomes
\begin{equation}\label{derexpconvexe2}
\frac{d}{dt} \int_{\rr^d} e^{\alpha \, \vert x \vert^2} \, d\mu_t(x) \leq
G \int_{\rr^d} (R^2 - \vert x \vert^2) \, e^{\alpha \, \vert x \vert^2} \, d\mu_t(x).
\end{equation} 

Let $p > 1$. The formula
\begin{eqnarray*}
\int_{\vert x \vert > p \, R} (R^2 - \vert x \vert^2) \, e^{\alpha \, \vert x \vert^2} \, d\mu_t(x)
& \leq &
R^2 ( 1- p^2) \int_{\vert x \vert > p \, R} e^{\alpha \, \vert x \vert^2} \, d\mu_t(x) \\
& = & 
R^2 ( 1- p^2) \left[ \int_{\rr^d} e^{\alpha \, \vert x \vert^2} \, d\mu_t(x) - 
\int_{\vert x \vert \leq p \, R} e^{\alpha \, \vert x \vert^2} \, d\mu_t(x) \right]
\end{eqnarray*}
leads to
$$
\int_{\rr^d} (R^2 - \vert x \vert^2) \, e^{\alpha \, \vert x \vert^2} \, d\mu_t(x)
\leq
\int_{\vert x \vert \leq p \, R} (R^2 p^2 - \vert x \vert^2) \, e^{\alpha \, \vert x \vert^2} \, d\mu_t(x)
+ R^2 (1-p^2) M_\alpha.
$$
by decomposing the integral on the sets $\{\vert x \vert \leq p \, R \}$ and
$\{\vert x \vert > p \, R \}$. From~\eqref{derexpconvexe2} we deduce
$$
(M_{\alpha})'(t) + \omega_1 \, M_\alpha(t) \leq \omega_2
$$
where $\omega_1$ and $\omega_2$ are positive constants. It follows that
$M_\alpha(t)$ remains bounded on $\rr_+$ if $M_\alpha(0)<+\infty$, and this concludes the
argument.
\end{proof}

\subsection{Time-regularity}

Now we study the time-regularity of $\mu_t$. 

\begin{proposition}\label{propcontW1-edp}
With the conventions of Subsection~\ref{subnotation}, 
for any $T<+\infty$ there exists a constant $C(T)$ such that
\begin{equation}\label{stCT} \forall s,t\in [0,T],\qquad
W_1(\mu_t,\mu_s) \, \leq \, C(T) \, \vert t-s \vert^{1/2}.
\end{equation}
\end{proposition}

\begin{remark} The exponent $1/2$ is natural in small time if no regularity assumption is made
on $\mu_0$; it can be improved if $t,s$ are assumed to be bounded below by some $t_0>0$.
Also, in view of the results of convergence to equilibrium recalled later on,
the constant $C(T)$ might be chosen independent of $T$ if 
$\beta >0, \beta+ 2\, \gamma >0$.
% or if the center of mass is fixed and $\beta +\gamma >0$.
\end{remark}

\begin{remark} A stochastic proof of~\eqref{stCT} is possible, via the study of continuity estimates
for $Y_t$, which in any case will be useful later on. But here we prefer to present
an analytical proof, to stress the fact that estimates in this section are purely
analytical statements.
\end{remark}

\begin{proof} 
Let $L$ be the linear operator $-\Delta - \nabla \cdot (\cdot \nabla V+\nabla(W\ast\mu_t))$,
and let $e^{-tL}$ be the associated semigroup: from our assumptions and estimates it
follows that it is well-defined, at least for initial data which admit a finite square
exponential moment. Of course $\mu_t=e^{-tL}\mu_0$. It follows that
\begin{align*} W_2 (\mu_s,\mu_t) = W_2(\mu_s, e^{-(t-s)L}\mu_s) & =
W_2 \left ( \int_{\rr^d} \delta_y \,d\mu_s(y),\: \int_{\rr^d} e^{-(t-s)L} \delta_y\,d\mu_s(y)
\right ) \\
& \leq \int_{\rr^d} W_2 (\delta_y, e^{-(t-s)L}\delta_y)\,d\mu_s(y). 
\end{align*}
Our goal is to bound this by $O(\sqrt{t-s})$. In view of Proposition~\ref{propT_1edp},
it is sufficient to prove that for all $a>0$,
\[ W_2^2 (\delta_y, e^{-(t-s)L}\delta_y) = O(t-s)\, O(e^{a|y|^2}). \]

This estimate is rather easy, since the left-hand side is just the variance
of the solution of a linear diffusion equation, starting with a Dirac mass at $y$ as
initial datum. Without loss of generality, we assume $s=0$, and write
$\tilde{\mu}_t:=e^{-tL}\delta_y$. For simplicity we write the computations in a
sketchy way, but they are not hard to justify. 

Since the initial datum is $\delta_y$, its square exponential moment $\tilde{M}_\alpha$
of order $\alpha$ is $e^{\alpha |y|^2}$. With an argument similar to the proof of 
Proposition~\ref{propT_1edp}(ii), one can show that 
\[ 0\leq t\leq T \Longrightarrow
\int e^{\alpha|x|^2} \,d\tilde{\mu}_t(x) \leq C(T) (1+\tilde{M}_\alpha) \leq
C(T)\, e^{\alpha|y|^2}. \]

Now, since
$|\nabla V|(x)=O(e^{a|x|^2})$, $a<\alpha$,
$|\nabla W\ast \mu_t|$ grows at most polynomially, and $\tilde{\mu}_t$ admits a
square exponential moment of order $\alpha$, we easily obtain
\[ \frac{d}{dt} \int x\,d\tilde{\mu}_t = - \int \nabla (V+W\ast\mu_t)\,d\tilde{\mu}_t
= \int O(e^{a|x|^2})\,d\tilde{\mu}_t = O(e^{\alpha |y|^2}); \]
\[ \frac{d}{dt} \int \frac{|x|^2}2\,d\tilde{\mu}_t = 
d - \int x\cdot \nabla (V+W\ast\mu_t)\,d\tilde{\mu}_t = O(e^{\alpha|y|^2}). \]
From these estimates we deduce that the time-derivative of the variance
$V(\tilde{\mu}_t):=\int |x|^2\,d\tilde{\mu}_t - (\int x\,d\tilde{\mu}_t)^2$ 
is bounded by $O(e^{b|y|^2})$ for any $b>0$. Since $\tilde{\mu}_0$ has zero variance, 
it follows that the variance of $\tilde{\mu}_t$ is $O(t e^{b|y|^2})$, which was our goal.
\end{proof}

\subsection{Regularity in phase space}

Regularity estimates will be useful for Theorem~\ref{thmreconstr}.
Equation~\eqref{edp} is a (weakly nonlinear) parabolic equation, for which
regularization effects can be studied by standard tools. Some limits to the
strength of the regularization are imposed by the regularity of $V$. So as not
to be bothered by these nonessential considerations, we shall assume strong
regularity conditions on $V$ here.
Then in Appendix \ref{sectionregulariteedp} we shall prove the following estimates:

\begin{proposition}\label{propregul}
With the conventions of Subsection~\ref{subnotation}, assume in addition that
$V$ has all its derivatives growing at most polynomially at infinity. Then,
for each $k\geq 0$ and for all $t_0>0$, $T>t_0$ there is a finite constant $C(t_0,T)$,
only depending on $t_0,T,k$ and a square exponential moment of the initial measure $\mu_0$,
such that the density $f_t$ of $\mu_t$ is of class ${\cal C}^k$, with
\[ \sup_{t_0\leq t\leq T} \|f_t\|_{{\cal C}^k} \leq C (t_0,T). \]
If moreover $\beta > 0$, $\beta+\gamma> 0$, 
then $C(t_0,T)$ can be chosen to be independent of $T$ for any fixed $t_0$.
\end{proposition}
%or if the center of mass is fixed and $\beta + \gamma > 0$, 

\begin{remark}
For regular initial data and under some adequate assumptions on $V$ and $W$, some
regularity estimates on $f_t/f_\infty$, where $f_\infty$ is the limit density
in large time, are established in~\cite[Lemma~6.7]{CMV04}. These estimates allow
a much more precise uniform decay, but are limited to just one derivative.
Here there will be no need for them.
\end{remark}

\subsection{Asymptotic behavior}

In the ``uniformly convex'' case when $\beta+\gamma>0$, 
the measure $\mu_t$ converges to a definite limit $\mu_\infty$ as $t\to\infty$. 
This was investigated in~\cite{Mal01,CMV03,CMV04}. The following statement
is a simple variant of~\cite[Theorems~2.1 and 5.1]{CMV03}. 

\begin{proposition} \label{propasympt}
With the conventions of Subsection~\ref{subnotation}, assuming that
$\beta >0, \beta+2 \, \gamma>0$, there exists a probability measure $\mu_\infty$ such that
\[W_2(\mu_t,\mu_\infty) \leq C e^{-\lambda t}, \qquad \lambda>0.\]
Here the constants $C$ and $\lambda$ only depend on the initial datum $\mu_0$.
\end{proposition}
%or that the center of mass is kept fixed by the equation and that
%$\beta + \gamma >0$, 

\section{The limit empirical measure}

Consider the random time-dependent measure
\begin{equation}\label{limem} 
\hat{\nu}^N_t := \frac1N \sum_{i=1}^N \delta_{Y_t^i}, 
\end{equation}
where $(Y_t^i)_{t\geq 0}$, $1\leq i\leq N$, are $N$ independent processes
solving the same stochastic differential equation
\[ dY_t^i = \sqrt{2}\,dB_t^i - [\nabla (V+W\ast \mu_t) ](Y_t^i)\,dt, \]
and such that the law of $Y_0^i$ is $\mu_0$. As we already mentioned,
for each $t$ and $i$, $Y_t^i$ is distributed according to the law $\mu_t$.
We call $\hat{\nu}^N_t$ the ``limit empirical measure'' because it is
expected to be a rather accurate description, in some well-chosen sense, 
of the empirical measure $\hat{\mu}^N_t$ as $N\to\infty$.

Our estimates on $\mu_t$, and the fact that $\hat{\nu}^N_t$ is the empirical
measure for \emph{independent} processes, are sufficient to imply good properties
of concentration of $\hat{\nu}^N_t$ around its mean $\mu_t$, as $N\to\infty$,
for each $t$. But later on we shall use some estimates about the time-dependent
measure (even to obtain a result of concentration for $\hat{\mu}^N_t$ with fixed $t$). 
To get such results, we shall study the time-regularity of $\hat{\nu}^N_t$.
Our final goal in this section is the following

\begin{proposition} \label{proppgd-eds}
With the conventions of Subsection~\ref{subnotation}, for any $T\geq 0$ there
are constants $C=C(T)$ and $a = a(T) >0$ such that the limit empirical measure~\eqref{limem}
satisfies
$$ \forall\Delta\in [0,T], \forall \eps>0,\qquad
\pp \Bigl[ \sup_{t_0\leq s, t \leq t_0+\Delta} W_1(\hat{\nu}_s^N, \hat{\nu}_t^N) > \eps\Bigr] 
\leq  \exp \left(-N(a \, \eps ^2 - C \, \Delta) \right).
$$
\end{proposition}

To prove Proposition~\ref{proppgd-eds}, we shall use a bit of classical stochastic calculus tools.

\subsection{SDE estimates}\label{sdeest}

In this subsection we establish the following estimates of time regularity for the
stochastic process $Y_t$: For all $T>0$, there exist positive constants $a$ and $C$ such that, for all 
$s,t,t_0, \Delta \in [0,T]$,

\[ \begin{array}{ccc}
\text{(i)} & \ee \, \vert Y_t - Y_s \vert^2 \leq C \vert t-s \vert \\ \\
\text{(ii)} & \ee \, \vert Y_t - Y_s \vert^4 \leq C \vert t-s \vert^2 \\ \\
\text{(iii)} & \ee \Bigl[ \ds \sup_{t_0 \leq s\leq t \leq t_0+\Delta} 
\exp \bigl(a \vert Y_t - Y_s \vert^2\bigr) \Bigr]\leq 1 + C \, \Delta.
\end{array} \]

\begin{proof}
We start with (i). We use It\^o's formula to write a stochastic equation on
the process $(\vert Y_t - Y_s \vert^2)_{t \geq s}$:
$$
\vert Y_t - Y_s \vert^2 = M_{s,t} + 2 \, d\,(t-s) - 
2 \int_s^t (\nabla V (Y_u) + \nabla W \ast\mu_u(Y_u) )\cdot (Y_u - Y_s) \, du,
$$
where $M_{s,t}$, viewed as a process depending on $t$, is a martingale with zero expectation. Hence
\begin{equation}\label{espcarre}
\ee \, \vert Y_t - Y_s \vert^2 =  
2 \, d\, (t-s) - 2 \int_s^t \ee \, (\nabla V (Y_u) + \nabla W \ast\mu_u(Y_u) )\cdot (Y_u - Y_s) \, du.
\end{equation}

On one hand
\begin{multline}\label{espcarre2}
\ee \Bigl\vert (\nabla V (Y_u) + \nabla W \ast\mu_u(Y_u) )\cdot (Y_u - Y_s) \Bigr\vert^2 
\leq
4 \Bigl( \ee \vert \nabla V (Y_u) \vert^2 + \ee \vert \nabla W \ast\mu_u(Y_u) \vert^2 \Bigr) \\ 
\Bigl(\ee \vert Y_u \vert^2 + \ee \vert Y_s \vert^2\Bigr).
\end{multline}
On the other hand, by Proposition~\ref{propT_1edp}, $\mu_u$ has a finite square exponential 
moment, uniformly bounded for $u\in [0,T]$. More precisely, 
there exist $\alpha >0$ and $M<+\infty$ such that 
$\ds \int e^{\alpha \vert x \vert^2} \, d\mu_u(x) \leq M$ for all $u \leq T$. 
Since by assumption $\vert \nabla W(z) \vert \leq L \, \vert z \vert$ and 
$\vert \nabla V(x) \vert = O(e^{\alpha \vert x \vert^2})$, 
we deduce 
$$ 
\sup_{s\leq u\leq T} \Bigl( \ee \, (\nabla V (Y_u) + \nabla W \ast\mu_u(Y_u) )\cdot (Y_u - Y_s) 
\Bigr) <+\infty.
$$ 
In view of~\eqref{espcarre}, it follows that there exists a constant $C=C(T)$ such that
$$
\ee \, \vert Y_t - Y_s \vert^2 \leq (2 d + C) \, (t-s).
$$
This concludes the proof of (i).

\medskip

To establish (ii), we perform a very similar computation.
For given $s$, let $Z_{s,t}:=(\vert Y_t - Y_s \vert^4)_{t \geq s}$.
Another application of It\^o's formula yields
\begin{multline*}
\ee \, Z_{s,t} = 4(2+d) \int_s^t \ee \vert Y_u - Y_s \vert^2 \, du
\\ - 4 \int_s^t \ee \vert Y_u - Y_s \vert^2 (Y_u - Y_s) \cdot 
\bigl(\nabla V(Y_u) + \nabla W\ast \mu_u(Y_u)\bigr) \, du.
\end{multline*}

On one hand, from (i),
$$
\int_s^t \ee \, \vert Y_u - Y_s \vert^2 \, du \leq 2 \, C \int_s^t (u-s) \, ds = C (t-s)^2.
$$
On the other hand
\begin{multline}
\int_s^t \ee \, \vert Y_u - Y_s \vert^2 (Y_u - Y_s) \cdot 
\bigl(\nabla V(Y_u) + \nabla W\ast \mu_u(Y_u)\bigr) \, du \\
\leq 
\left( \int_s^t \ee \, Z_{s,u} \, du \right)^{3/4} 
\left( \int_s^t \ee \, \bigl\vert \nabla V(Y_u) + 
\nabla W\ast \mu_u(Y_u) \bigr\vert^4 \, du \right)^{1/4}\
\end{multline}
by H\"older's inequality. But again, since the measures $\mu_t$ admit a bounded 
square exponential moment, $\ee \, \vert \nabla V(Y_u) + \nabla W\ast \mu_u(Y_u) \vert^4$
is bounded on $[0,T]$. We conclude that
\begin{equation}\label{espquadr}
\ee \, Z_{s,t} \leq C \left( (t-s)^2 + (t-s)^{1/4} 
\left( \int_s^t \ee \, Z_{s,u} \, du \right)^{3/4} \right).
\end{equation}
Then, with $C$ standing again for various constants which are independent of $s$ and $t$,
$$
\ee \, Z_{s,u} \leq C (\ee \vert Y_u \vert^4 + \ee \vert Y_s \vert^4 ) 
\leq 2 \, C \sup_{0 \leq t \leq T} \int \vert x \vert^4 \, d\mu_u(x) \leq C;
$$
so, from~\eqref{espquadr},
$$
\ee \, Z_{s,t} \leq C((t-s)^2 + (t-s)^{1/4} (t-s)^{3/4}) \leq C (t-s),
$$
and by~\eqref{espquadr} again we successively  obtain
$$
\ee \, Z_{s,t} \leq C(t-s)^{7/4},
$$
and finally
$$
\ee \, Z_{s,t} \leq C(t-s)^2.
$$
This concludes the proof of (ii).

\medskip

We finally turn to the proof of (iii). Without real loss of generality,
we set $t_0=0$. We shall proceed as in 
the proof of Proposition~\ref{propT_1edp}, and prove the existence of some constant $C$ 
and some continuous positive function $a$ on $\rr_+$ such that
\begin{equation}\label{espexp}
\ee \left(\sup_{0\leq s \leq t \leq \Delta\leq T} \exp \left( a(t) \vert Y_t - Y_s \vert^2 \right)
\right) \leq 1 + C \, \Delta.
\end{equation}

Let $a(t)$ be a smooth function, and
\[ Z_{s,t}  := e^{ a(t) \vert Y_t - Y_s \vert^2}.\]
By It\^o's formula,
\begin{multline*}
Z_{s,t} = 1 + M_{s,t} \\ + \int_s^t \left[ 2 a(u) 
\Bigl( d + 2 a \vert Y_u - Y_s \vert^2 - (\nabla V+ \nabla W \ast \mu_u )(Y_u)
\cdot (Y_u - Y_s) \Bigr) + a'(u) \vert Y_u - Y_s \vert^2 \right] \, Z_{s,u} \, du
\end{multline*}
where
$$
M_{s,t} := \int_s^t a(u) \, (Y_u - Y_s) \, Z_u \, dB_u.
$$
For each $s$, $M_{s,t}$, viewed as a stochastic process in $t$, is a martingale.

By Young's inequality, for any $b>0$,
$$
-2 \bigl(\nabla V+ \nabla W \ast \mu_u\bigr)(Y_u) \cdot (Y_u - Y_s) 
\leq b \bigl\vert Y_u - Y_s \bigr\vert^2 + \frac{1}{b}
\bigl\vert \nabla V+ \nabla W \ast \mu_u \bigr\vert^2(Y_u).
$$
So, by letting
$$
A_u :=  a(u) \left[ 2\, d + \frac{1}{b} \bigl\vert \nabla V(Y_u) + 
\nabla W \ast \mu_u (Y_u) \bigr\vert^2 \right]
$$
and
$$
B(u) := a'(u) + 4 \, a^2(u) + b \, a(u)
$$
we obtain
$$
Z_{s,t} \leq 1 + M_{s,t} + \int_s^t [A_u + B(u) \vert Y_u - Y_s \vert^2] \, Z_{s,u} \, du.
$$

We choose $a$ in such a way that the function $B$ is identically zero, that is
$$
a(u) = e^{-b\, u} \left(\frac{1}{a(0)} + 4 \, \frac{1 - e^{-b\, u}}{b} \right)^{-1},
$$
where $a(0)$ is to be fixed later. Then
$$
Z_{s,t} \leq 1 + M_{s,t} + \int_s^t A_u \, Z_{s,u} \, du
$$
from which it is clear that
\begin{equation}\label{ineqZst}
\ee \sup_{s \leq t \leq \Delta} Z_{s,t} \leq 
1 + \ee \sup_{s \leq t \leq \Delta} M_{s,t} + \int_s^\Delta \ee \, A_u \, Z_{s,u} \, du.
\end{equation}

By Cauchy-Schwarz and Doob's inequalities,
\begin{equation}\label{doob}
\left( \ee \sup_{s \leq t \leq \Delta} M_{s,t} \right)^2 
\leq \ee \, \left\vert \sup_{s \leq t \leq \Delta} M_{s,t} \right\vert^2 
\leq 2 \, \sup_{s \leq t \leq \Delta} \ee \, \vert M_{s,t} \vert^2.
\end{equation}
Also, by It\^o's formula and the Cauchy-Schwarz inequality again,
\begin{multline}\label{doob2}
\ee \, \vert M_{s,t} \vert^2  = \int_s^t a(u)^2 \, \ee \vert Y_u - Y_s \vert^2 \, Z_{s,u}^2 \, du 
\\
\leq  \frac12 \int_s^t a(u)^2 \, \left( \ee\,  \vert Y_u - Y_s \vert^4 \right)^{1/2}
\left( \ee \, Z_{s,u}^4 \right)^{1/2} \, du.
\end{multline}
In view of (ii), there exists a constant $C$ such that
\begin{equation}\label{espquadr2}
\ee \, \vert Y_u - Y_s \vert^4 \leq C \, (u-s)^2.
\end{equation}
Furthermore,
\begin{equation}\label{Z4}
\ee \, Z_{s,u}^4 = \ee \exp 4 \, a(u) \vert Y_u - Y_s \vert^2 \leq
\left( \ee  \exp 16 \, a(u) \vert Y_u\vert^2 \right)^{1/2} 
\left( \ee  \exp 16 \, a(u) \vert Y_s \vert^2 \right)^{1/2}.
\end{equation}
Recall from Proposition~\ref{propT_1edp} that there exist constants $M$ and $\alpha >0$ such that
$$
\sup_{s \leq u \leq \Delta} \int e^{\alpha \vert y \vert ^2} \, d\mu_u(y) \leq M.
$$
If we choose $a(0) \leq \alpha/16$, the decreasing property
of $a$ will ensure that $a(u)\leq \alpha/16$ for all $u\in [0,\Delta]$, and
$$
\ee  \exp 16 \, a(u) \vert Y_u\vert^2  
\left( = \int e^{16 \, a(u) \vert y \vert^2} \, d\mu_u(y) \right) \leq M.
$$
Then, from~\eqref{Z4},
$$
\sup_{s \leq u \leq \Delta} \ee \, Z_{s,u}^4 \, \leq M.
$$
Now, from~\eqref{doob2} and~\eqref{espquadr2} we deduce
$$
\sup_{s \leq t \leq \Delta} \ee \, \vert M_{s,t} \vert^2 \leq C \, (t-s)^2.
$$
Combining this with \eqref{doob}, we conclude that
$$
\ee \sup_{s \leq t \leq \Delta} M_{s,t} \leq C \, \Delta.
$$

In the same way, we can prove that $\ee \, ( A_t  Z_{s,t})$ is bounded for $t\in[s,\Delta]$ 
by bounding $\ee \, Z_{s,t}^2$ and $\ee \, A_t^2$.  
This concludes the proof of~\eqref{espexp}, and therefore of (iii) above.
\end{proof}

\subsection{Time-regularity of the limit empirical measure}

We are now ready to prove Proposition~\ref{proppgd-eds}.

On one hand
$$
W_1(\hat{\nu}^N_s, \hat{\nu}^N_t) \leq \frac{1}{N} 
\sum_{i=1}^{N} \vert Y_t^i - Y_s^i \vert,
$$
so
\begin{equation}\label{Vshi}
\pp \, \left[ \sup_{0 \leq s \leq t \leq \Delta} W_1(\hat{\nu}_s^N, \hat{\nu}_t^N) > \eps\right] 
\leq  \pp \left[ \frac1{N} \sum_{i=1}^{N} V^i > \eps \right]
\end{equation}
where
$$
V^i := \sup_{0 \leq s \leq t \leq \Delta} \vert Y_t^i - Y_s^i \vert.
$$

By Chebyshev's exponential inequality and the independence of the $(Y_t^i-Y^i_s)$,
$$
\pp \left[ \frac1{N} \sum_{i=1}^{N} V^i > \eps \right] \leq
\exp\left( -N \sup_{\zeta \geq 0} \, \bigl[ \eps \zeta - \ln \ee \exp ( \zeta V^1) \bigr] \right).
$$
But, for any given $\zeta$ and $\omega \geq 0$,
$$
\ee \exp ( \zeta V^1) \leq \ee \exp 
\left(\zeta \bigl(\frac{\omega^2 + (V^1)^2}{2 \, \omega}\bigr)\right)
\leq \exp \frac{\zeta \, \omega}{2} \, \ee \exp \frac{\zeta}{2 \, \omega} (V^1)^2.
$$
Let $\ds \omega = \frac{\zeta}{2 \, a}$, so that $\ds \frac{\zeta}{2 \, \omega} = a$.
Then, from estimate (iii) in Subsection~\ref{sdeest},
$$
\ee \exp \frac{\zeta}{2 \, \omega} (V^1)^2 \leq 1 + C \, \Delta,
$$
uniformly in $s$ and $\Delta$. Hence, for any $\zeta>0$,
$$
\ee \exp ( \zeta V^1) \leq \ee \exp \frac{\zeta^2}{4 \, a} \, (1 + C \, \Delta).
$$

Consequently,
\begin{eqnarray*}
\pp \left[ \frac{1}{N} \sum_{i=1}^{N} V^i > \eps \right] 
& \leq &
\exp \left(-N \sup_{\zeta\geq 0} \, 
\bigl[ \, \eps \zeta - \frac{\zeta^2}{4} - \ln (1 + C \, \Delta) \bigr] \right) \\
& = & \exp \Bigl(-N \, [ a \, \eps^2 - \ln (1 + C \, \Delta) ] \Bigr) \\
& \leq & \exp \left(-N \, [ a \, \eps^2 - C \, \Delta] \right).
\end{eqnarray*}
The proof of Proposition~\ref{proppgd-eds} follows by~\eqref{Vshi}.

\section{Coupling}

We now (as is classical) reduce the proof of convergence for $\hat{\mu}^N_t$
to a proof of convergence for the empirical measure $\hat{\nu}^N_t$ constructed
on the auxiliary independent system $(Y_t^i)$. The final goal of this section
is the following estimate.

\begin{proposition} \label{propcouplage}
With the conventions of Subsection~\ref{subnotation},
\[ W_1(\hat{\mu}^N_t, \mu_t) \leq \Gamma \int_0^t e^{-\alpha (t-s)} W_1(\hat{\nu}^N_s,
\mu_s)\,ds + W_1(\hat{\nu}^N_t,\mu_t), \]
where $\Gamma$ is defined by~\eqref{Gamma}, and $\alpha:= \beta+ 2\min(\gamma,0)$.
\end{proposition}

\begin{proof}
For the sake of simplicity we give a slightly sketchy proof. We couple the stochastic
systems $(X_t^i)$ and $(Y_t^i)$ by assuming that (i) $X_0^i=Y_0^i$ and
(ii) both systems are driven by the \emph{same} Brownian processes $B_t^i$.
In particular, for each $i\in \{1,\ldots,N\}$,
the process $X_t^i - Y_t^i$ satisfies the equation
\begin{equation}\label{eqX-Y}
d(X_t^i - Y_t^i) = - \bigl(\nabla V(X_t^i) - \nabla V(Y_t^i)\bigr) \, dt - 
\Bigl(\nabla W \ast \hat{\mu}_t^N (X_t^i) - \nabla W\ast \mu_t(Y_t^i)\Bigr) \, dt.
\end{equation}

From~\eqref{eqX-Y} we deduce
\begin{multline} \label{startoption}
\frac{1}{2} \frac{d}{dt} \vert X_t^i - Y_t^i \vert^2 =
- \, (\nabla V(X_t^i) - \nabla V(Y_t^i)) \cdot (X_t^i - Y_t^i) \\ 
- \, \bigl(\nabla W \ast \hat{\mu}_t^N (X_t^i) - 
\nabla W\ast \mu_t(Y_t^i)\bigr) \, \cdot\,  (X_t^i - Y_t^i). 
\end{multline}
Our convexity assumption on $V$ implies
\[ - (\nabla V(X_t^i) - \nabla V(Y_t^i)) \cdot (X_t^i - Y_t^i) 
\leq - \beta \vert X^i_t - Y^i_t \vert^2;
\]
so the main issue consists in the treatment of the quantity
$\nabla W \ast \hat{\mu}_t^N (X_t^i) - \nabla W\ast \mu_t(Y_t^i)$
appearing in the right-hand side of~\eqref{startoption}. 
There are (at least) two options here. The first one consists in writing
\begin{multline}\label{opt1}
\nabla W \ast \hat{\mu}_t^N (X_t^i) - \nabla W\ast \mu_t(Y_t^i)= \\
(\nabla W\ast\hat{\mu}_t^N-\nabla W\ast\mu_t)(X_t^i)
+ \bigl(\nabla W \ast \mu_t (X_t^i) - \nabla W\ast \mu_t(Y_t^i)\bigr);
\end{multline}
while the second one consists in forcing the introduction of $\hat{\nu}_t^N$ 
as follows:
\begin{multline}\label{opt2} 
\nabla W \ast \hat{\mu}_t^N (X_t^i) - \nabla W\ast \mu_t(Y_t^i)=\\
\frac1{N}\sum_{j=1}^{N} \bigl[\nabla W(X_t^i-X_t^j) - \nabla W(Y_t^i-Y_t^j)\bigr]
-(\nabla W\ast \hat{\nu}^N_t-\nabla W\ast\mu_t)(Y^i_t).
\end{multline}
Both options are interesting and lead to slightly different computations.
Since both lines of computations might be useful in other contexts, 
we shall sketch them one after the other. The second option leads to
better bounds, but at the price of more complications (in particular,
we shall need to sum over the index $i$ at an early stage).
\med

\noindent{\bf First option}:
We start as in~\eqref{opt1}.
In view of our assumption on $D^2W$, the Lipschitz norm of $\nabla W(X_t^i - \, \cdot \,)$ 
is bounded by $\Gamma$. Therefore, by the Kantorovich-Rubinstein dual formulation~\eqref{KR},
$$
\Bigl| \nabla W \ast (\hat{\mu}_t^N -\mu_t)(X_t^i) \Bigr| = 
\left\vert \int_{\rr^d} \nabla W(X_t^i - y) \, d(\hat{\mu}_t^N - \mu_t)(y) \right\vert \leq 
\Gamma \, W_1(\hat{\mu}_t^N,\mu_t),
$$
and then our assumptions on $V$ and $W$ imply
$$
\frac{1}{2} \frac{d}{dt} \vert X_t^i - Y_t^i \vert^2  
\leq - (\gamma + \beta) \, \vert X_t^i - Y_t^i \vert^2
+ \Gamma \, W_1(\hat{\mu}_t^N,\mu_t) \, \vert X_t^i - Y_t^i \vert.
$$

In other words, $\vert X_t^i - Y_t^i \vert$ satisfies the differential inequality
$$
\frac{d}{dt} \vert X_t^i - Y_t^i \vert + (\beta +\gamma) \,  \vert X_t^i - Y_t^i \vert
\leq \Gamma \,   W_1(\hat{\mu}_t^N,\mu_t)
$$
($X_t^i$ and $Y_t^i$ separately are not Lipschitz functions of $t$, but their
difference is). Hence, by Gronwall's lemma,
$$
\vert X_t^i - Y_t^i \vert \leq \Gamma \int_{0}^{t} e^{-(\beta+\gamma)(t-s)} \, 
W_1(\hat{\mu}_s^N,\mu_s) \, ds.
$$

Now we sum over $i$; by convexity of the distance $W_1$ and 
triangular inequality, we obtain
\begin{eqnarray*}
W_1(\hat{\mu}_t^N,\hat{\nu}_t^N)  
& \leq &
\frac{1}{N} \sum_{i=1}^{N} \vert X_t^i - Y_t^i \vert 
\leq  \Gamma \int_{0}^{t} e^{-(\beta+\gamma)(t-s)} \, W_1(\hat{\mu}_s^N,\mu_s) \, ds \\
& \leq & \Gamma \int_{0}^{t} e^{-(\beta+\gamma)(t-s)} \, \bigl[W_1(\hat{\mu}_s^N,\hat{\nu}_s^N) 
+ W_1(\hat{\nu}_s^N,\mu_s) \bigr] \, ds.
\end{eqnarray*}

By using Gronwall's lemma again, we deduce
$$
W_1(\hat{\mu}_t^N,\hat{\nu}_t^N)  \leq  
\Gamma \int_{0}^{t} e^{-(\beta+\gamma -\Gamma) (t-s)} \, W_1(\hat{\nu}_s^N,\mu_s) \, ds.
$$
By applying the triangular inequality for $W_1$, we conclude to
the validity of Proposition~\eqref{propcouplage}, only with $\alpha$ 
replaced by the (a priori smaller) quantity $\beta+\gamma-\Gamma$.
\med

\noindent{\bf Second option:} Now we start with~\eqref{opt2}.
This time we sum over $i$ right from the beginning:
\[
\frac{1}{2} \frac{d}{dt} \sum_{i=1}^{N} \vert X_t^i - Y_t^i \vert^2 =
- \sum_{i=1}^{N} (\nabla V(X_t^i) - \nabla V(Y_t^i)) \cdot (X_t^i - Y_t^i) 
-\frac{1}{N} \sum_{i,j=1}^{N} (A^{ij}_t + B^{ij}_t)
\]
where
\[
A^{ij}_t
= (\nabla W(X^i_t-X^j_t) - \nabla W(Y^i_t - Y^j_t)) \, \cdot\,  (X_t^i - Y_t^i) 
\]
and
\[
B^{ij}_t = (W(Y^i_t - Y^j_t) - \nabla W\ast \mu_t(Y_t^i)) \, \cdot\,  (X_t^i - Y_t^i).
\]

Since $\nabla W$ is an odd function and $D^2 W(x) \geq \gamma I$ for all $x \in \rr^d$, we have
\begin{multline*}
A^{ij}_t + A^{ji}_t
= \bigl(\nabla W(X^i_t-X^j_t) - \nabla W(Y^i_t-Y^j_t)\bigr) \, 
\cdot\,  \bigl((X_t^i - X_t^j) - (Y_t^i-Y^j_t)\bigr) \\
\geq \gamma \bigl\vert (X_t^i - X_t^j) - (Y_t^i-Y^j_t) \bigr\vert^2,
\end{multline*}
whence
\[
- \sum_{i,j=1}^{N} A^{ij}_t \leq - \frac{\gamma}{2} \sum_{i,j=1}^{N} 
\bigl\vert (X_t^i - X_t^j) - (Y_t^i-Y^j_t) \bigr\vert^2 
\leq - 2 N \gamma^{-} \sum_{i=1}^{N} \vert X^i_t - Y^i_t \vert^2
\]
where $\gamma^{-} = \min(\gamma,0)$.

Then
\[
- \sum_{j=1}^{N} B^{ij}_t = 
-(X^i_t - Y^i_t) \cdot (\nabla W \ast \hat{\nu}^N_t (Y^i_t) - \nabla W \ast \mu_t(Y^i_t)).
\] 
Our assumption on $D^2W$ implies that the Lipschitz norm of 
$\nabla W(Y_t^i - \, \cdot \,)$ is bounded by $\Gamma$; so, 
by the Kantorovich-Rubinstein dual formulation~\eqref{KR},
$$
\Bigl| \nabla W \ast (\hat{\nu}_t^N -\mu_t)(Y_t^i) \Bigr| = 
\left\vert \int_{\rr^d} \nabla W(Y_t^i - y) \, d(\hat{\nu}_t^N - \mu_t)(y) \right\vert \leq 
\Gamma \, W_1(\hat{\nu}_t^N,\mu_t).
$$

Collecting all terms we finally obtain
$$
\frac{1}{2} \frac{d}{dt} \sum_{i=1}^N \vert X_t^i - Y_t^i \vert^2  
\leq - (\beta + 2 \gamma^{-}) \sum_{i=1}^N  \vert X_t^i - Y_t^i \vert^2
+ \Gamma \sum_{i=1}^N \vert X_t^i - Y_t^i \vert  W_1(\hat{\nu}_t^N,\mu_t).
$$

Then, since 
$\ds \sum_{i=1}^N \vert X_t^i - Y_t^i \vert \leq \left( N \sum_{i=1}^N \vert X^i_t - Y^i_t \vert^2\right)^{1/2}$,
the function $\ds y(t) := \left(\frac{1}{N} \sum_{i=1}^N  \vert X^i_t - Y^i_t \vert^2 \right)^{1/2}$ 
satisfies the differential inequality
\[
y'(t) + (\beta + 2 \gamma^{-}) y(t) \leq \Gamma \, W_1(\hat{\nu}^N_t, \mu_t),
\]
so that
\[
\left( \frac{1}{N} \sum_{i=1}^N  \vert X^i_t - Y^i_t \vert^2 \right)^{1/2}
\leq  \Gamma \int_{0}^{t} e^{-(\beta+2 \gamma^{-})(t-s)} \, W_1(\hat{\nu}_s^N,\mu_s)  \, ds.
\]

The conclusion follows by triangular inequality again since
\[
W_1(\hat{\mu}_t^N,\hat{\nu}_t^N)  \leq  
W_2(\hat{\mu}_t^N,\hat{\nu}_t^N) \leq \left( \frac{1}{N} \sum_{i=1}^N  \vert X^i_t - Y^i_t \vert^2 \right)^{1/2}.
\]
\end{proof}

\begin{remark} Not only does the ``second option'' in the proof lead to better
bounds, it also provides an estimate of the distance between $\hat{\mu}$
and $\hat{\nu}$ in the $W_2$ distance, which is stronger
than the $W_1$ distance. However, we do not take any advantage of
this refinement.
\end{remark}

\section{Conclusion}

In this section, we paste together all the estimates established in the previous
sections, so as to prove Theorems~\ref{thmconc-eds} to~\ref{thmreconstr}.

\subsection{Concentration estimates}

We start with the proof of Theorem~\ref{thmconc-eds}. 
By $C$ we shall denote various constants depending on $T$, on our assumptions
on $V$ and $W$, and also on $\int e^{\alpha|x|^2}\,d\mu_0(x)$, for some
$\alpha>0$.

From Proposition~\ref{propcouplage},
\[ \sup_{0\leq t\leq T} W_1(\hat{\mu}^N_t, \mu_t) \leq 
(\Gamma e^{|\alpha|T}+1) \sup_{0\leq t\leq T} W_1(\hat{\nu}^N_t, \mu_t)\,ds.\]
In particular, there is a constant $C$ such that
\begin{equation}\label{mutonuhat}
\pp \, \left[\sup_{0 \leq t \leq T} W_1(\hat{\mu}_t^N, \mu_t) > \eps\right] \leq
\pp \, \left[\sup_{0 \leq t \leq T} W_1(\hat{\nu}_t^N, \mu_t) > \tilde{\eps} \right],
\qquad \tilde{\eps}=\frac{\eps}{C}.
\end{equation}

From Corollary~\ref{sqexpT1} and Theorem~\ref{thmconc-Tp2} we know that
$$
\sup_{0\leq t\leq T} \pp \, [ W_1(\hat{\nu}_t^N, \mu_t) > \tilde{\eps}] 
\leq e^{-K \, N \, \tilde{\eps}^2}
$$
for all $t \in [0,T]$, $N \geq  N_0\,\max(\tilde{\eps}^{-(d'+2)},1)$ ($d'>d$).
The issue now is to ``exchange'' $\sup$ and $\pp$ in this estimate.
As we shall see, this is authorized by the continuity estimates on $\hat{\nu}_t^N$ and
$\mu_t$.

Let $\Delta >0$ (to be fixed later on), and let $M$ be the integer part of $T/\Delta \, + \,1$.
We decompose the interval $[0,T]$ as
\[ [0,T] = [0,\Delta] \cup [\Delta,2\Delta] \cup \ldots \cup [(M-1)\Delta, T]
\subset \bigcup_{h=0}^{M-1} [h\Delta, (h+1)\Delta]. \]
Proposition~\ref{propcontW1-edp} guarantees that, if $\Delta\leq a\tilde{\eps}^2$ 
for some $a$ small enough, then
\begin{equation}\label{holder}
h\Delta\leq t\leq (h+1)\Delta \Longrightarrow 
W_1(\mu_t,\mu_{h\Delta}) \leq \frac{\tilde{\eps}}{2}.
\end{equation}

Then, by triangular inequality and~\eqref{holder},
$$
\ds \pp \, \left[\sup_{0 \leq t \leq T} W_1(\hat{\nu}_t^N, \mu_t) > \tilde{\eps} \right] 
$$
$$
\leq \ds \pp \,  \left[\sup_{h=0, \dots, M-1} \: \sup_{h\Delta \leq t \leq (h+1)\Delta} 
W_1(\hat{\nu}_t^N, \mu_t) > \tilde{\eps} \right] 
$$
\begin{multline*}
\leq \pp \,  \left[\sup_{h=0, \dots, M-1}\: \sup_{h\Delta \leq t \leq (h+1)\Delta}  
W_1(\hat{\nu}_t^N, \nu_{h\Delta}^N) + 
\sup_{h=0, \dots, M-1} W_1(\hat{\nu}_{h\Delta}^N, \mu_{h\Delta}) \right.
\\ \left. + \sup_{h=0, \dots, M-1}\: \sup_{h\Delta \leq t \leq (h+1)\Delta} 
W_1(\mu_{h\Delta}, \mu_t)> \tilde{\eps}) \right]
\end{multline*}
$$
\leq \pp \,  \left[\sup_{h=0, \dots, M-1}\: \sup_{h\Delta \leq t \leq (h+1)\Delta}  
W_1(\hat{\nu}_t^N, \hat{\nu}_{h\Delta}^N) +
\sup_{h=0, \dots, M-1} W_1(\hat{\nu}_{h\Delta}^N, \mu_{h\Delta}) 
> \frac{\tilde{\eps}}{2} \right],
$$
which can be bounded by
$$
\ds \pp \,  \left[\sup_{h=0, \dots, M-1}\: \sup_{h\Delta \leq t \leq (h+1)\Delta}  
W_1(\hat{\nu}_t^N, \hat{\nu}_{h\Delta}^N) > \frac{\tilde{\eps}}{4} \right]
\, + \, \pp \, \left[\sup_{h=0, \dots, M-1} 
W_1(\hat{\nu}_{h\Delta}^N, \mu_{h\Delta}) > \frac{\tilde{\eps}}{4} \right].
$$

By Corollary~\ref{sqexpT1} and Theorem~\ref{thmconc-Tp2}, there exist some constants $C$ and $N_0$ such that
$$
\pp \, \left[ W_1(\hat{\nu}_{h\Delta}^N, \mu_{h\Delta}) \geq 
\frac{\tilde{\eps}}{4} \right] \leq \exp( - C \, N \, \tilde{\eps}^2)
$$
for all $h=0, \dots, M-1$, and $N \geq N_0 \max(\tilde{\eps}^{-(d'+2)},1)$. Hence 
\begin{equation}\label{1stterm}
\ds \pp \,  \left[\sup_{h=0, \dots, M-1} W_1(\hat{\nu}_{h\Delta}^N, \mu_{h\Delta}) 
> \frac{\tilde{\eps}}{4} \right]
\leq \sum_{h=0}^{M-1} \pp \, \left[ W_1(\hat{\nu}_{h\Delta}^N, \mu_{h\Delta}) 
> \frac{\tilde{\eps}}{4} \right] \leq M \exp( - C \, N \, \tilde{\eps}^2).
\end{equation}

On the other hand, from Proposition~\ref{proppgd-eds} we deduce
$$
\pp \, \left[\sup_{h\Delta \leq t \leq (h+1)\Delta}  
W_1(\hat{\nu}_t^N, \hat{\nu}_{h\Delta}^N) > \frac{\tilde{\eps}}{4} \right]
\leq \exp \left( -N \bigl( \frac{a}{4} \, \tilde{\eps}^2 - C \, \Delta\bigr) \right)
$$
for all $h=0, \dots, M-1$ and $\tilde{\eps} > 0$, so
\begin{equation}\label{presqueest}
\pp \, \left[\sup_{h=0, \dots, M-1}\: \sup_{h\Delta \leq t \leq (h+1)\Delta}  
W_1(\hat{\nu}_t^N, \hat{\nu}_{h\Delta}^N) > \frac{\tilde{\eps}}{4}
\right] \leq M \exp \left( -N \bigl( \frac{a}{4} \, \tilde{\eps}^2 - C \, \Delta\bigr) \right).
\end{equation}
We can assume that $\Delta \leq \frac{a}{8 \, C} \tilde{\eps}^2$, and $M\le CT/\tilde{\eps}^2 +1$;
then we can bound the right-hand side of~\eqref{presqueest} by
\begin{equation}\label{2ndterm}
M \exp \left( -\frac{a}{8} \, N \,\tilde{\eps}^2 \right)
\leq C \left(1+\frac{T}{\tilde{\var}^2}\right)\; 
\exp \left( -\frac{a}{8} \, N \,\tilde{\eps}^2 \right)
\end{equation}

From~\eqref{1stterm} and~\eqref{2ndterm} we deduce that, for $\Delta$  small enough
(depending on $\var$!), 
\begin{equation}\label{2terms}
\pp \, \left[\sup_{0 \leq t \leq T} W_1(\hat{\nu}_t^N, \mu_t) > \tilde{\eps} \right] \leq
2 \, C \left(1+\frac{T}{\tilde{\var}^2}\right) \,  \exp( - K \, N \, \tilde{\eps}^2)
\end{equation}
for $N \geq N_0 \max(\tilde{\eps}^{-(d'+2)},1)$. 
So we deduce from~\eqref{2terms} that
$$
\pp \, \left[\sup_{0 \leq t \leq T} W_1(\hat{\nu}_t^N, \mu_t) > \tilde{\eps}\right] 
\leq \exp \left( \ln\left( C \left( \frac{T}{\tilde{\eps}^2}+1
\right) \right) - K \, N \, \tilde{\eps}^2 \right)
,
$$
where again $C, K$ stand for various positive constants, and 
$N\geq \max(N_0\,\eps^{-(d'+2)},1)$.
This concludes the proof of Theorem~\ref{thmconc-eds}.

\subsection{Uniform in time estimates}

Now, we shall focus on the case when $\beta >0, \beta + 2 \, \gamma >0$ is positive, and derive
Theorem~\ref{thmconc-eds-unif} by a slightly refined estimate.

Let us start again from the bound
\[ W_1(\hat{\mu}^N_t, \mu_t) \leq \Gamma \int_0^t e^{-\alpha (t-s)}
W_1(\hat{\nu}^N_s, \mu_s)\,ds + W_1(\hat{\nu}^N_t, \mu_t) \]
where $\alpha := \beta + 2 \min(\gamma, 0)$ is positive.
Let $\Delta>0$ (to be fixed later on), and $k$ be the integer part
of $t/\Delta$. If $W_1(\hat{\mu}^N_t,\mu_t)$ is larger than $\var$, then
\[ \begin{cases} \text{either } \ds W_1(\hat{\nu}^N_t,\mu_t) \geq \frac{\var}2\\
\\ \text{or  } \ds \exists j\in \{0,\ldots,k\}; \quad
\int_{j\Delta}^{(j+1)\Delta} e^{-\alpha(t-s)} W_1(\hat{\nu}^N_s, \mu_s)\,ds
\geq \frac{\var}{2^{k+2-j}\Gamma}. \end{cases} \]
Indeed, $(\var/2)+\sum_{j\leq k} (\var/2^{k+2-j})\leq \var$.
As a consequence,
\[ \begin{cases} \text{either } \ds W_1(\hat{\nu}^N_t,\mu_t) > \frac{\var}2\\
\\ \text{or  } \ds \exists j\in \{0,\ldots,k\}; \quad
\sup_{j\Delta\leq s\leq (j+1)\Delta} W_1(\hat{\nu}^N_s, \mu_s)
> \frac{\var\alpha  e^{\alpha [t-(j+1)\Delta]}}{2^{k+2-j}\Gamma}. \end{cases} \]
Since, for $t\in [j\Delta, (j+1)\Delta]$,
\[ \frac{e^{\alpha [t-(j+1)\Delta]}}{2^{k+2-j}} \ge 
\frac{e^{\alpha (k-j-1)\Delta}}{2^{k-j+2}} =
\left ( \frac1{4e^{\alpha \Delta}} \right ) \left 
(\frac{e^{\alpha\Delta}}2 \right)^{k-j},\]
we conclude to the existence of a constant $C$ such that
\begin{multline} \label{estjj+1}
\pp \left[ W_1(\hat{\mu}^N_t, \mu_t) > \var\right] \leq
\pp \left[ W_1(\hat{\nu}^N_t,\mu_t) > \frac{\var}2 \right]
\\ + \sum_{j=0}^k \pp \left[\sup_{j\Delta\leq s\leq (j+1)\Delta} 
W_1(\hat{\nu}^N_s, \mu_s) > C \eps \left(\frac{e^{\alpha\Delta}}2 \right)^{k-j} \right].
\end{multline}

We already know that the first term in the right-hand side in~\eqref{estjj+1}
is bounded by $e^{-\lambda N \var^2}$ for some constant $\lambda>0$,
and so we focus on the other terms.

In the proof of Theorem~\ref{thmconc-eds}, we have established that
there are constant $C$ and $\lambda$, depending on $\Delta$ and on bounds
on square exponential moments for $\mu_0$, such that
\begin{equation}\label{Clambdabds} \pp \left[\sup_{0\leq s\leq \Delta} 
W_1(\hat{\nu}^N_s, \mu_s)> \delta\right] \leq C \left(1+\frac{\Delta}{\delta^2}\right)\: e^{-\lambda N\delta^2}.
\end{equation}
Proposition~\ref{propT_1edp} guarantees that these square exponential
bounds also hold true for $\mu_t$, uniformly in $t$. Thus we can
apply~\eqref{Clambdabds} with $\mu_{j\Delta}$ taken as initial datum,
and get
\begin{equation}\label{Clambdabdsj} \pp \left[\sup_{j\Delta\leq s\leq (j+1)\Delta} 
W_1(\hat{\nu}^N_s, \mu_s)> \delta\right] \leq C e^{-\lambda N\delta^2},
\end{equation}
as soon as $N\geq N_0\max(\delta^{-(d'+2)},1)$.
\med

We now use~\eqref{Clambdabdsj} to bound the sum appearing in the 
right-hand side of~\eqref{estjj+1}.
Choose $\Delta$ large enough that
\[ \theta:= \frac{e^{\alpha\Delta}}2 > 1. \]
Applying~\eqref{Clambdabdsj} with $\delta$ replaced by
$C\theta^{j-k}\var$, 
we can bound the sum in the right-hand side of~\eqref{estjj+1} by
\[ C \sum_{j=0}^k \exp\bigl(-K \theta^{2(k-j)} N\var^2\bigr)\]
for $N \geq N_0 \max(\eps^{-(d'+2)}, 1)$, where $C$, $K$  and $N_0$ are again positive constants. 
Since again $\theta$ is larger than~1, there is a constant $a>0$ such that $\theta^{2(k-j)}\geq a(k-j)$,
so the sum above is bounded by
\[  C \left (e^{-KN\var^2} + \sum_{\ell=1}^\infty e^{-K\ell N\var^2}\right )
\leq C \left (e^{-KN\var^2} + \frac{e^{-KN\var^2}}{1-e^{-KN\var^2}} \right).\]
If $N_0$ is large enough, our assumption $N\geq N_0\max(\var^{-(d'+2)},1)$ 
implies that $e^{-KN\var^2}$ is always less than $1/2$, so that
the above sum can be bounded by just $C e^{-KN\var^2}$. This concludes
the proof of the first point of Theorem~\ref{thmconc-eds-unif}.

\bigskip

The second point is proved by writing 
\begin{eqnarray*}
W_1(\hat\mu^N_t,\mu_\infty)&\le& W_1(\hat\mu^N_t,\mu_t)+W_1(\mu_t,\mu_\infty)\\
&\le& W_1(\hat\mu^N_t,\mu_t)+Ce^{-\lambda t}
\end{eqnarray*}
successively by the triangular inequality for Wasserstein distance and use of Proposition \ref{propasympt}. 
Then the result follows  from the uniform estimate obtained above.

\subsection{Data reconstruction}

We finally consider Theorem~\ref{thmreconstr}.
Proposition~\ref{propregul} ensures that, as $t\to\infty$, 
$f_t$ is uniformly bounded in $C^k$, where $k$ is arbitrarily large.
Since $f_t$ converges to $f_\infty$ as $t\to\infty$, we deduce that
$f_\infty$ is Lipschitz. Then Theorem~\ref{thmconc-eds-unif}
and Proposition~\ref{propreconstr} together imply Theorem~\ref{thmreconstr}.

\appendix

\section{Metric entropy of a probability space}\label{sectionentropiemetrique}

We now prove the covering result used in Section~\ref{sectionpreuve}, 
as a particular case of a more general estimate. Let $E$ be a Polish space,
we look for an upper bound on the number $\NN_p(E,\delta) := m(\PP(E),\delta)$ of
balls of radius $\delta$ in Wasserstein distance $W_p$ needed to cover the
space $\PP(E)$ of probability measures on $E$. We use the same strategy as
in~\cite[Exercise~6.2.19]{DZ98}, where the L\'evy distance is used instead of
the Wasserstein distance.

\bthm \label{Th-balls} 
Let $(E,d)$ be a Polish space with finite diameter $D$. For any $r>0$,
define $N(E,r)$ as the minimal number of balls needed to cover $E$ by
balls of radius $r$. Then there exists a numerical constant $C$ 
such that for all $p\geq 1$ and $\delta \in (0, D)$, the space
$\PP(E)$ can be covered by $\NN_p(E,\delta)$ balls of radius 
$\delta$ in $W_p$ distance, with
\begin{equation} \label{nballs}
\NN_p(E,\delta)\le \left( \frac{C D}{\delta}\right)^{pN\left(E,\frac{\delta}{2}\right)}.
\end{equation}  
\nthm

\begin{remark}
The $W_p$ distance between any two probability measures on $E$ is at most $D$, so, 
for all $\delta\geq D$, we have the trivial estimate $\NN_p(E,\delta) = 1$.
\end{remark}

\bprf

Let $r >0$, and let $\{x_j\}_{1\leq j\leq N(E,r)}$ be such that $E$ is covered by the 
balls $B(x_j,r)$ with centers $x_j \in E$ and radius $r$. For simplicity we
shall write $N=N(E,r)$.
\med

In a {\bf first step} we prove that for any $\mu\in\PP(E)$ there exist nonnegative real numbers
$(\beta_j)_{1 \leq j \leq N}$, with $ \ds \sum_{j=1}^{N} \beta_j =1$, such that
$$
W_p(\mu, \tilde \mu) \leq r, \qquad \tilde \mu := \sum_{j=1}^{N} \beta_j \delta_{x_j}.
$$

For this we first replace the balls $B(x_j,r)$'s by the sets $\tilde B_j$'s defined by
$$
\forall j, \qquad \tilde B_j =  B(x_j,r) \setminus \bigcup_{k \leq j-1} B(x_k,r),
$$
so that $E$ is partitioned into the $\tilde B_j$'s. 
Next define
$$
\beta_j = \mu[\tilde B_j].
$$
It is easy to check that the required properties are fulfilled. Indeed, we may
transport $\mu$ onto $\ds \tilde \mu = \sum_{j=1}^{N} \beta_j \delta_{x_j}$ 
by sending all $x$'s in $\tilde B_j$ onto $x_j$, for each $j=1, \dots, N$: the cost of
this transport is bounded by
$\sum_{j=1}^{N} r^p \mu(\tilde B_j)  = r^p$.

\medskip

In the {\bf second step} we introduce an integer $K$ (whose value will be made more precise 
later on), and consider the set
$$
\CC_K := \left\{ \sum_{j=1}^{N} \alpha_j \delta_{x_j}; \quad
(\alpha_j)_{1 \leq j \leq N} \in \AA_K \right\} \subset \PP(E),
$$
where $\AA_K$ is the set of all $N$-tuples 
$(\alpha_j)_{1\leq j\leq N}$, such that each $\alpha_j$ is of the form $k_j/K$, $k_j\in\N$, and
$\ds \sum_{j=1}^{N} \alpha_j  =1$.

Given a probability measure $\ds \tilde \mu = \sum_{i=1}^{N} \beta_i \delta_{x_i}$
(where $(\beta_i)_i$ does not necessarily belong to $\AA_K$), there exists
$\mu'$ in $\CC_K$ such that
\begin{equation}\label{mumu'}
W_p(\mu', \tilde \mu) \leq D \left( \frac{N}{K} \right)^{1/p}.
\end{equation}

To prove~\eqref{mumu'}, we define $n_j$ as the integer part 
$[K \beta_j]$ of $K \beta_j$ and $J$ as the first integer such that
$$
\sum_{j=1}^{J} (n_j+1) + \sum_{j =J+1}^{N} n_j = K.
$$
Since $\ds \sum_{j=1}^{N} \beta_j =1$, it is clear that $J\leq N$.
Then we define a measure $\mu'\in \CC_K$ by
$\ds \mu' = \sum_{j=1}^{N} \alpha_j \delta_{x_j}$, where
$$
\alpha_j  = \left\{ \begin{array}{ll}
                           \frac{n_j+1}{K} & \mbox{for $ j=1, \dots, J$}\\
                           \frac{n_j}{K} & \mbox{for $ j=J+1, \dots, N.$}   
                          \end{array}
                  \right.
$$

Let us bound the distance between $\mu$ and $\mu'$. For that
we gradually define a transport plan between $\tilde \mu$ and $\mu'$ in the following way: 
first of all, at each point $x_i$, the mass $n_i/K$ stays in place.
Then, the remaining masses $\beta_i-n_i/K$ are redistributed as follows:
all the remaining mass at $x_1,\ldots, x_\ell$ is brought to $x_1$,
together with possibly a bit of mass at $x_{\ell+1}$, until a total mass
$1/K$ has been added at location $x_1$ (for $\ell$ large enough).
If $J\geq 2$, then we again bring mass from $x_{\ell+1},\ldots$, until
another mass $1/K$ has been added at $x_2$. We carry on until all the mass
at $x_J$ has been used, thus building a transport plan 
$(\pi_{ij})_{1 \leq i,j \leq N}$ which sends $\tilde \mu$ onto $\mu'$, 
in such a way that $\ds \pi_{ii} \geq \frac{n_i}{K}$ for all $i$. Hence,
$$
\sum_{j \neq i} \pi_{ij} \leq \beta_i - \pi_{ii} = \beta_i - \frac{n_i}{K} \leq \frac{1}{K},
$$ 
and this plan yields an upper bound on the Wasserstein distance:
$$
W_p^p(\tilde \mu, \mu')
\leq \sum_{i,j=1}^{N} d(x_i, x_j)^p \, \pi_{ij} 
=  \sum_{i=1}^{N} \, \sum_{j \neq i}  d(x_i, x_j)^p \, \pi_{ij}
\leq  N \frac{D^p}{K}.
$$

\medskip

To summarize the first two steps: for any $\mu$ in $\PP(E)$ there exists $\mu' \in \CC_K$ such that
$$
W_p(\mu, \mu') \leq \ r + D \left( \frac{N}{K} \right)^{1/p}.
$$
In other words, the family $\Bigl(B \bigl(\mu', r+ D(N/K)^{1/p}\bigr)\Bigr)_{\mu' \in \CC_K}$ 
covers $\PP(E)$.

\bigskip

In the {\bf third step} we choose some suitable $K$ and $r$ for a given $\delta$.

We first choose $K$ in such a way that $r$ and $D (N/K)^{1/p} $ 
 have the same order of magnitude, for instance
$$
K = \left[ N\left(\frac{D}{r}\right)^p \right] + 1.
$$

Then
$$
r+D\,(N/K)^{1/p} \leq 2 r ,
$$
and the balls $B \bigl(\mu', r+ D(N/K)^{1/p}\bigr)$ have radius at most $\delta$ if
$$
r = \frac{\delta}{2}.
$$

\medskip

Now $K$ and $r$ are fixed, $N=N(E,\delta/2)$, and we just have to estimate the cardinality 
$\sharp\CC_K$ of $\CC_K$. For this we first note that
$$
\sharp \CC_K = \frac{(K+N-1)!}{(K-1)!N!} = \frac{(K+N-1) \dots K}{N!}
$$
Without loss of generality, we have assumed $\delta<D$, so $K>N$. Then
$K< \dots < K+N-1 < 2K$, and hence
$$
\sharp \CC_K \leq \frac{(2K)^N}{N!} \leq \left(\frac{2Ke}{N} \right)^N.
$$

Since $N \geq 1$ and $2D \geq \delta$, we can write
$$
K \leq N \left(\frac{2D}{\delta}\right)^p + 1 \leq 2N \left(\frac{2D}{\delta}\right)^p,
$$
and we deduce
$$
\sharp \CC_K \leq \left(C\,\frac{D}{\delta} \right)^{pN\bigl(E,\frac{\delta}{2}\bigr)}
$$
with $C = 2(4e)^{1/p} \leq 8 e$.

\medskip

Consequently, we have covered $\PP(E)$ by the 
$\ds \left(C \frac{D}{\delta} \right)^{pN\bigl(E,\frac{\delta}{2}\bigr)}$ 
balls $(B(\mu',\delta))_{\mu' \in \CC_K}$ with radius $\delta$.
This concludes the argument.

\nprf

\bigskip

In the particular case when $E$ is the Euclidean ball $B_R$ of radius $R$ in $\rr^d$, we have
\begin{equation}\label{nballsR} 
N(B_R,r)\le k\left(\frac{R}{r}\right)^d
\end{equation} 
for some constant $k$. To see this, one may for instance consider the balls with center 
in the lattice $\frac{r}{\sqrt{d}} \zz^d$ in $\rr^d$. Then Theorem~\ref{Th-balls} yields the bound
$$
\NN_p(B_R,\delta) \le \left( C\,\frac{R}{\delta}\right)^{pk\left(\frac{R}{\delta} \right)^d},
$$
which is used in the present paper.

\section{Regularity estimates on the limit PDE}\label{sectionregulariteedp}

In this appendix we study solutions to the limit equation
\begin{equation}\label{edp2}
\partial_t \rho = \Delta \rho + \nabla. (\rho (V + W \ast \rho)), \qquad t \geq 0, \quad x \in \rr^d
\end{equation}
and establish the regularity results stated in Proposition \ref{propregul}. 
Following the method in~\cite{DV00}, we shall measure the regularity
in terms of $L^2$-Sobolev spaces
\[
H^s(\rr^d) = \Bigl\{ u \in L^2(\rr^d); \partial^{\alpha} u
\in L^2(\rr^d), \, \alpha \in \nn^d, \, \vert \alpha \vert \leq s \Bigr\}
\qquad (s\in\nn).
\]
Our main result is as follows.

\bthm\label{thmregul}
Let $V$ and $W$ such that all their partial derivatives 
$\ds \partial^{\alpha} V$ and $\ds \partial^{\alpha} W$ are continuous 
and grow at most polynomially at infinity, for any multi-index 
$\alpha \in \nn^d$ with $\vert \alpha \vert \leq s+1$.
Let $a, E>0$ and let $\rho_0$ be a probability density such that
\[
\int_{\rr^d} e^{a \vert x \vert^2}\, d\rho_{0}(x) \leq E.
\] 
Then, there exists a continuous function $f:(0,+\infty)\to(0,+\infty)$,
only depending on $d$, $s$, $V$, $W$, $a$ and $E$,
such that any classical solution $\rho = \rho(t,x)$ to \eqref{edp2},
starting from $\rho_{0}$, satisfies
\[
\bigl\Vert \rho(t,\cdot) \bigr\Vert_{H^s(\rr^d)} \leq f(t).
\]
\nthm

\begin{proof}

For the sake of simplicity we only give a formal proof, which can be turned rigorous by means of regularization
arguments.

Let then $\rho = (\rho(t,.))_{t \geq 0}$ be a solution of
\[
\partial_t \rho = \Delta \rho + \nabla. \left(\rho (V + W \ast \rho)\right), \qquad t \geq 0, \quad x\in \rr^d;
\]
we rewrite the equation as
\[
\partial_t \rho = \sum_{i=1}^{d} \partial_{ii} \rho +
\partial_i \left[ \rho \partial_i \phi\right],
\]
where $\partial_i = \partial^{e_i}$ if $e_i$ is the $i$-th vector of 
the canonical base of $\rr^d$, and
\[
\phi(t,x) = V(x) + W\ast \rho(t,x).
\]

Let $\alpha \in \nn^d$ be given. By integration by parts and Cauchy-Schwarz 
inequality,
\begin{multline*}
\frac{1}{2} \frac{d}{dt}\int_{\rr^d} \left\vert \partial^{\alpha} \rho \right\vert^2
= \int_{\rr^d} \partial^{\alpha} \rho \; \partial_t \left( \partial^{\alpha} \rho \right) = 
\int_{\rr^d} \partial^{\alpha} \rho \; \partial^{\alpha} \left(\partial_t \rho\right) 
\\
= \sum_{i=1}^{d} \int_{\rr^d} \partial^{\alpha} \rho  \;  \partial^{\alpha}
\left( \partial_{ii} \rho \, + \, \partial_i \left[ \rho \, \partial_i \phi \right]  \right)
\\
= - \sum_{i=1}^d \int_{\rr^d} \left\vert \partial^{\alpha+e_i} \rho \right\vert^2 
\; + \; \int \partial^{\alpha+e_i} \rho  \; \partial^{\alpha} \left[ \rho \partial_i \phi\right] 
\\
\leq - \sum_{i=1}^d \int_{\rr^d}\left\vert \partial^{\alpha +e_i} \rho \right\vert^2 
\; + \; \sum_{i}^d \left[ \int_{\rr^d} \left\vert \partial^{\alpha+e_i} \rho \right\vert^2 \right]^{1/2} \; 
\left[ \sum_{\beta \leq \alpha} C_{\alpha, \beta} \int_{\rr^d} \left\vert \partial^{\alpha -\beta+e_i} \phi
\; \partial^{\beta} \rho\right\vert^2 \right]^{1/2}
\\
\leq -\frac{1}{2} \sum_{i=1}^d \int_{\rr^d}\left\vert \partial^{\alpha +e_i} \rho \right\vert^2 
\; + \; \sum_{\beta \leq \alpha} C_{\alpha, \beta} \int_{\rr^d} \left\vert \partial^{\alpha -\beta+e_i} \phi
\; \partial^{\beta} \rho\right\vert^2.
\end{multline*}

By summing over $\alpha \in \nn^d$ with $\vert \alpha \vert \leq s$, we find
\[
\frac{d}{dt} \sum_{\vert \alpha \vert \leq s} \int_{\rr^d} \left\vert \partial^{\alpha} \rho \right\vert^2
\leq - \sum_{\vert \alpha \vert \leq s} \sum_{i=1}^d \int_{\rr^d}\left\vert \partial^{\alpha +e_i} \rho \right\vert^2 
\; + \; \sum_{\vert \alpha \vert \leq s} \sum_{\beta \leq \alpha} C_{\alpha, \beta} \int_{\rr^d} 
\left\vert \partial^{\alpha -\beta+e_i} \phi \; \partial^{\beta} \rho\right\vert^2.
\]

Given $T >0$, by Proposition \ref{propT_1edp} there exist constants 
${\hat a}$ and ${\hat E}$, depending only on $d$, $a$, $E$ and $T$, such that
\begin{equation}\label{momentexp}
\int e^{{\hat a} \vert x\vert^2} d\rho(t,x) \leq {\hat E}
\end{equation}
for all $t \in [0,T]$. In particular, it follows from our assumptions
on the derivatives of $V$ and $W$ that all 
$\ds \left\vert \partial^{\alpha -\beta+e_i} \phi\right\vert^2$
terms are bounded by some polynomial in $|x|$, uniformly in $t \in [0,T]$.

Let $\langle x\rangle:=\sqrt{1+|x|^2}$. For $k, s\geq 0$, we introduce the
weighted norms
\[
\Vert u \Vert_{H_k^s} := \left( \sum_{\vert \alpha \vert \leq s} 
\int_{\rr^d} \langle x\rangle^k 
\, \vert \partial^{\alpha} u(x) \vert^2 \, dx \right)^{1/2}
\]
and
\[ \Vert u\Vert_{L^1_k} := \int_{\rr^d} \langle x\rangle^k\, \vert u(x)\vert
\,dx.\]
Then for any $s \in \nn$ and $T \geq 0$ there exist $k$ and $C \geq 0$ 
such that
\begin{equation}\label{ineq}
0\leq t\leq T\Longrightarrow\qquad
\frac{d}{dt} \Vert u \Vert_{H^s}^2 \leq -\Vert u \Vert_{H^{s+1}}^2 + C \, \Vert u \Vert_{H_k^s}^2.
\end{equation} 

We shall prove later on the following interpolation lemma:

\begin{lemma}\label{lemmeinterpol}
Given $d \geq 1$, $s\in \nn$ an $k \geq 0$, there exist nonnegative 
constants $C(d,s,k)$ and $h(d,s,k)$, and $\theta(d,s) \in (0,1)$ 
such that for all $u \in L^{1}_{\infty}(\rr^d) \cap H^{s+1}(\rr^d)$,
\[
\Vert u \Vert_{H^s_k} \leq 
C(d,s,k) \Vert u \Vert_{L^1_{h(d,s,k)}}^{1 - \theta(d,s)} \, \Vert u \Vert_{H^{s+1}}^{\theta(d,s)}.
\]
\end{lemma}

\medskip

Then, again from \eqref{momentexp}, all $\Vert u \Vert_{L^1_{h(d,s,k)}}(t)$ 
norms are bounded on $[0,T]$, so from \eqref{ineq} and Lemma~\ref{lemmeinterpol}
 there exists some constants $C$ such that
\[
\frac{d}{dt} \Vert u \Vert_{H^s}^2 \leq -\Vert u \Vert_{H^{s+1}}^2 + C \, \Vert u \Vert_{H^{s+1}}^{2 \, \theta}
\leq - \frac{1}{2} \Vert u \Vert_{H^{s+1}}^2 + C 
\leq - C  \Vert u \Vert_{H^s}^{2/\theta} + C.
\]

In other words $\ds A(t) = \Vert u \Vert_{H^s}^2(t)$ satisfies on $[0,T]$ the
differential inequality
\begin{equation}\label{nash}
A'(t) + c \, A(t)^{p} \leq C
\end{equation}
for some constants $c, C \geq 0$ and $p = 1/{\theta} >1$ depending only on 
$d$, $a$, $E$, $s$ and $T$.

\medskip

Let us distinguish two cases.
If $A(0) \leq 1$, then we only use the inequality $A'(t) \leq C$
to make sure that
\[
A(t) \leq A(0) + C t \leq 1 + C T
\]
for any $t \in [0,T]$.

If on the other hand $A(0) \geq 1$, we deduce from \eqref{nash} that
\[
A'(t) + c \, A(t)^{p} \leq C A(t),
\]
as long as $A(t) \geq 1$, so that $D(t) := A(t)^{1-p}$ satisfies the inequality
\[
D'(t) + (p-1) \, C D(t) \geq (p-1) \, c
\]
which integrates to
\[
D(t) \geq D(0) e^{(1-p)Ct} + \frac{c}{C} (1 - e^{(p-1)C t}) \geq \frac{c}{C} (1 - e^{(p-1)t}).
\]
As a consequence, as long as $A(t) \geq 1$, we have
\[
A(t) \leq (c/C)^{1/1-p} (1 - e^{(p-1)t})^{1/(1-p)}.
\]

\medskip

In the end, we have obtained an a priori bound on 
$\ds A(t) = \int \left\vert \partial^{\alpha} \rho \right\vert^2 (t)$ 
for $t \in (0,T]$, depending only on $d, s, a, E$ and $T$, but not on 
the initial value $A(0)$. Then the proof can be concluded by an 
approximation argument.
\end{proof}

\bigskip

\begin{proof}[Proof of Lemma~\ref{lemmeinterpol}]
We proceed by induction on $s$.

In the {\bf first step} we prove the result for $s=0$. Given $d \geq 1$ and $a \in (0,1]$, we write
\[
\int_{\rr^d}  \langle x\rangle^{k} \, \vert u(x) \vert^2 \, dx =
\int_{\rr^d} \langle x\rangle^{k} \, \vert u(x) \vert^a \, \vert u(x) \vert^{2-a}\, dx,
\]
so, by H\"older's inequality, 
\[
\Vert u \Vert_{L^2_k}^2 \leq \Vert u \Vert_{L^1_{\frac{k}{a}}}^a \; \Vert u \Vert_{L^{\frac{2-a}{1-a}}}^{2-a}
\]
(with $\ds \frac{2-a}{1-a} = \infty$ if $a=1$).
Then by Sobolev embedding,
\[
\Vert u \Vert_{L^2_k}^2 \leq C(d,a) \Vert u \Vert_{L^1_{\frac{k}{a}}}^a \; \Vert u \Vert_{H^1}^{2-a},
\]
where $a=1$ if $d=1$, $a$ is arbitrary in $(0,1)$ if $d=2$,
and $\ds a = \frac{4}{d+2}$ if $d \geq 3$, that is,
\[
\Vert u \Vert_{L^2_k} \leq C(d) \Vert u \Vert_{L^1_{\frac{k}{a}}}^{1 - \theta(d)} \; \Vert u \Vert_{H^1}^{\theta(d)}
\]
where $\ds \theta(1) = \frac{1}{2}$, 
any $\ds \theta(2) \in \left(\frac{1}{2}, 1 \right)$ for $d=2$, and
$\ds \theta(d) = \frac{d}{d+2}$ for $d \geq 3$.

\medskip

In the {\bf second step} we let $s \geq 1$ and assume by induction that there exist some constants $C(d,s-1,k),
h(d,s-1,k) \geq 0$ and $\theta(d,s-1) \in (0,1)$ such that for all $u\in L^1_{\infty}(\rr^d)\cap H^s(\rr^d)$:
\[
\Vert u \Vert_{H^{s-1}_k} \leq 
C(d,s-1,k) \Vert u \Vert_{L^1_{h(d,s-1,k)}}^{1 - \theta(d, s-1)} \, \Vert u \Vert_{H^s}^{\theta(d, s-1)}.
\]

Let then $u \in L^1_{\infty}(\rr^d) \cap H^{s+1}(\rr^d)$.

Given $\alpha \in \nn^d$ with $\vert \alpha \vert = j$ and $1 \leq j \leq s$, 
we split $\alpha$ into $\alpha = \alpha_1 + \alpha_2$ with 
$\vert \alpha_2 \vert =1$, and integrate by parts:
\begin{multline*}
\Vert \partial^{\alpha} u \Vert^2_{L^2_k} \leq  k \, \Vert \partial^{\alpha_1} u \Vert_{L^2_{2k-2}}
\; \Vert \partial^{\alpha} u \Vert_{L^2} + \Vert \partial^{\alpha_1} u \Vert_{L^2_{2 k}} 
\; \Vert \partial^{\alpha + \alpha_2} u \Vert_{L^2} \\
\leq (k +1) \Vert \partial^{\alpha_1} u \Vert_{L^2_{2 k}} \; \sup_{\vert \alpha \vert \leq j+1}
\Vert \partial^{\alpha} u \Vert_{L^2},
\end{multline*}
whence
\begin{multline*}
\sup_{\vert \alpha \vert = j} \Vert \partial^{\alpha} u \Vert^2_{L^2_k}
\leq (k+1)  \sup_{\vert \alpha \vert = j-1} \Vert \partial^{\alpha} u \Vert_{L^2_{2k}}
\sup_{\vert \alpha \vert \leq j+1} \Vert \partial^{\alpha} u \Vert_{L^2} \\
\leq (k+1)  \sup_{\vert \alpha \vert \leq s-1} \Vert \partial^{\alpha} u \Vert_{L^2_{2k}}
\sup_{\vert \alpha \vert \leq s+1} \Vert \partial^{\alpha} u \Vert_{L^2}.
\end{multline*}

Since this holds for any $ 1 \leq j \leq s$ we obtain
\[
\sup_{1 \leq \vert \alpha \vert \leq s} \Vert \partial^{\alpha} u \Vert^2_{L^2_k}
\leq (k+1)  \sup_{\vert \alpha \vert \leq s-1} \Vert \partial^{\alpha} u \Vert_{L^2_{2k}}
\sup_{\vert \alpha \vert \leq s+1} \Vert \partial^{\alpha} u \Vert_{L^2}.
\]

Moreover
\[
\Vert u \Vert^2_{L^2_k} \leq \Vert u \Vert_{L^2_{2\, k}} \; \Vert u \Vert_{L^2} 
\leq \sup_{\vert \alpha \vert \leq s-1} \Vert \partial^{\alpha} u \Vert_{L^2_{2 k}}
\sup_{\vert \alpha \vert \leq s+1} \Vert \partial^{\alpha} u \Vert_{L^2},
\]
so that finally
\[
\Vert u \Vert_{H^s_k}^2 \leq (k+1)\Vert u \Vert_{H^{s-1}_{2 k}} \; \Vert u \Vert_{H^{s+1}}.
\]

Then, by induction hypothesis, 
\[
\Vert u \Vert_{H^s_k}^2 \leq (k+1) \,
C(d,s-1,2k) \Vert u \Vert_{L^1_{h(d,s-1,2k)}}^{1 - \theta(d, s-1)} \, \Vert u \Vert_{H^s}^{\theta(d, s-1)}
\; \Vert u \Vert_{H^{s+1}},
\]
whence
\[
\Vert u \Vert_{H^s_k} \leq C(d,k,s) \Vert u \Vert_{L^1_{h(d,s,k)}}^{1 - \theta(d, s)} \,
\Vert u \Vert_{H^{s+1}}^{\theta(d,s)}
\]
where $\ds \theta(d,s) = \frac{1}{2 - \theta(d,s-1)} \in (0,1)$ and $h(d,s,k) = h(d,s-1,2k) \geq 0$.
This concludes the argument.

\end{proof}

\med
\noindent
{\bf Acknowledgments:} The authors thank M. Ledoux for his relevant comments and his interest during 
the preparation of this work, as well as providing Reference \cite{GZ86}.
\med

\nocite{*}
%\bibliography{./sanov}

\end{document}